\documentclass[12pt,fleqn]{amsart}

\usepackage[all]{xy}
\usepackage{tikz-cd}
\usepackage{amscd}
\usepackage{amsfonts}
\usepackage{color}
\usepackage{enumerate}
\usepackage{amssymb}
\newtheorem{theorem}{Theorem}[section]
\newtheorem{lemma}[theorem]{Lemma}
\newtheorem{proposition}[theorem]{Proposition}
\newtheorem{problem}[theorem]{Problem}
\newtheorem{corollary}[theorem]{Corollary}

\theoremstyle{definition}
\newtheorem{definition}[theorem]{Definition}

\newtheorem{remark}[theorem]{Remark}

\theoremstyle{remark}

\numberwithin{equation}{section}

\newcommand{\cA}{\mathcal A}

\newcommand{\cC}{\mathcal C}

\newcommand{\cD}{\mathcal D}

\newcommand{\cL}{\mathcal L}
\newcommand{\cN}{\mathcal N}
\newcommand{\cU}{\mathcal U}

\newcommand{\cZ}{\mathcal Z}

\newcommand{\fA}{\mathfrak A}
\newcommand{\fB}{\mathfrak B}
\newcommand{\fC}{\mathfrak C}

\newcommand{\cF}{\mathcal F}
\newcommand{\cH}{\mathcal H}
\newcommand{\ult}{\protect{\rm ult}}
\newcommand{\clop}{\protect{\rm clop}}
\newcommand{\la}{\langle}
\newcommand{\ra}{\rangle}
\newcommand{\sub}{\subseteq}
\newcommand{\eps}{\varepsilon}
\newcommand{\er}{\mathbb R}
\newcommand{\sm}{\setminus}
\newcommand{\con}{\mathfrak c}

\newcommand{\vf}{\varphi}
\newcommand{\oo}{{\omega_1}}

\newcommand{\btu}{\bigtriangleup}

\newcommand{\MA}{\mbox{\sf MA}}
\newcommand{\CH}{\mbox{\sf CH}}

\newcommand{\wh}{\widehat}
\newcommand{\pe}{\it p}

\newcommand{\mt}{{\rm mt}}
\newcommand{\bil}[2]{{\left\langle\kern0ex #1,#2
		\kern0ex\right\rangle}}

\newcommand{\ol}{\overline}
\newcommand{\stevo}{Todor\v{c}evi\'c}

\newcommand{\om}{\omega}
\newcommand{\zok}{[0,1]^\kappa}
\newcommand{\ttk}{2^\kappa}
\newcommand{\mk}{\lambda_\kappa}
\newcommand{\kr}{ \,{\bf\dot{}}\, }
\newcommand{\MS}{\hspace{0.1cm}\protect{\sc MS}}
\newcommand{\ct}{\tau_{\rm conv}}
\usepackage[geometry]{ifsym}
\DeclareMathOperator\BTU{\mbox{\BigTriangleUp}}
\begin{document}

\title[A survey on $P(K)$ spaces]{A survey on topological properties of $P(K)$ spaces}

\author[G.\ Plebanek]{Grzegorz Plebanek}
\address{Instytut Matematyczny\\ Uniwersytet Wroc\l awski\\ Pl.\ Grunwaldzki 2/4\\
50-384 Wroc\-\l aw\\ Poland} \email{grzegorz.plebanek@math.uni.wroc.pl}

\thanks{Partially  supported by the NCN (National Science Centre, Poland),
under the Weave-UNISONO call in the Weave programme 2021/03/Y/ST1/00.}

\begin{abstract}
Given a compact space $K$, we denote by $P(K)$ the space of all Radon
probability measures on $K$, equipped with the $weak^\ast$ topology inherited
from $C(K)^\ast$. For nonmetrizable compacta $K$ even  basic properties
of $P(K)$ spaces depend of additional axioms of set theory.
We discuss here  older and quite recent results on the subject. 
\end{abstract}

\date{}

\maketitle

\begin{flushright}
\begin{minipage}{23em}
{\em  To the memory of {\em Ja\'s}, Kamil Duszenko (1986-2014)\\
Ten years after}
\end{minipage}
\end{flushright}
\vspace{3ex}

\section{Introduction}
In the sequel, $K$ always stands for a compact Hausdorff space.
By $P(K)$ we denote the space of all probability Radon measures on $K$. Thus every 
$\mu\in P(K)$ is a probability measure defined on the Borel $\sigma$-algebra $Bor(K)$ which is inner-regular, i.e.\ 
\[ \mu(B) = \sup\{\mu(F): F\sub B, F \mbox{ compact}\},\]
for every $B\in Bor(K)$. 

Every $\mu\in P(K)$ acts by integration on $C(K)$, the Banach space of continuous functions, as a norm-one functional.
Since $P(K)\sub C(K)^\ast$,  the space $P(K)$ is naturally equipped with the $weak^\ast$ topology, the weakest one
making all the mappings
\[ P(K)\ni \mu \mapsto \int_K f\; {\rm d}\mu, \quad f\in C(K),\]
continuous. 
Then $P(K)$ is a compact space and its topological properties are the
main subject of this survey.

There are  many measure-theoretic issues that arise in the context of Banach spaces; Rodr\'{\i}guez \cite{Ro17}
offers a survey on such problems. However, here we wish  to confine our considerations to results 
directly related to topological structure of compacta of the form $P(K)$ and the main objective
is to discuss the interplay between properties of $K$ and related properties of $P(K)$. In turn,
 features of $P(K)$ often quite directly influence some properties of the Banach space $C(K)$ but those
are mentioned {\em en passant}. 
Let us also mention here that, from a more abstract perspective, $P$ is a covariant functor in the category of compact spaces, and
some properties of $P(K)$ spaces, such as cellularity (not discussed here), can be computed by categorical approach, see
Fedorchuk and \stevo\ \cite{FT97}.

If $K$ is a metrizable compact space then the topological structure of $P(K)$ is readily understood: for infinite $K$
the space $P(K)$ is  homeomorphic to the Hilbert cube $[0,1]^\omega$. For nonmetrizable compacta $K$, however,
even questions about basic topological properties of $P(K)$ become subtle and many natural problems
turn out to be independent of the usual axioms of set theory. 
Years ago Fedorchuk \cite{Fe91}  discussed a variety of such results in his extensive survey. 
However, substantial progress has occurred since then, and we aim to present these advancements here. 
In particular, there are  interesting,  very recent results that we discuss in section 7.

Typically, we offer no proofs, sometimes sketching basic ideas, but there are a few  exceptions.

\begin{enumerate}[(A)]
\item In section \ref{fccc} we reprove some nearly classical results on measures on Corson compacta, trying to convince the reader that designing the family of closed-and-open sets of the target compactum may be far quicker
than building it by means of inverse systems.
\item 
In section \ref{mtc} we present a measure-theoretic lemma devised  by Richard Haydon --- that auxiliary result
 serving as a shortcut to theorems regarding continuous surjections from
 $K$ or $P(K)$ onto Tikhonov cubes.
The lemma was already described and published,  or better say concealed, in \cite{Pl02a}, the paper
bearing a suspicious title. 
\item In section \ref{sp} we outline Fremlin's result on uniformly distributed sequences on the Cantor cube
$2^\con$. 
Although this theorem has been around for some time, it appears to have been solely published in Fremlin's treatise and remains relatively obscure, even among specialists,
see e.g.\   Mercourakis and Vassiliadis \cite{MV23}. 
\item In section \ref{groef} we present a result stating that for two non-scattered compact
spaces $K$ and $L$,  the fact that $C(K\times L)$ is not a Grothendieck space
can be attested  by a sequence of nonatomic measures from $P(K\times L)$.
\end{enumerate}

Let us mention here that our approach to  (D) has been recently developed in \cite{PRS} where we
prove, in particular, that $C[0,1]$ is isomorphic to a complemented subspace of the Banach space $C(\beta\omega\times\beta\omega)$.

We use very little from Banach Space Theory; let us  mention Diestel's book \cite{Di84} and Albiac and Kalton \cite{AK06}
as entirely sufficient sources. We use several, standard by now, additional axioms of set theory and we refer
to Jech \cite{Jech} in such cases. Let us remark, however, that many problems 
hinge on one's response to the following question: 

\begin{quote}
 {\em Given $\mu\in P(K)$ and an uncountable family of
Borel sets with positive measure,  is there a point $x\in K$ belonging to uncountably many of those sets?}
\end{quote}

 --- see subsection \ref{cm} for details.

We mention below a number of classical results from general topology and those can be found in Engelking's
monograph \cite{En}.
Concerning topological measure theory, we refer to Fremlin's treatise \cite{Fr3,Fr4,Fr5} (where
there is nearly everything) or to the more concise book by Bogatchev \cite{Bo07}. 

In fact, we take the opportunity and
offer a mini-course on finite Radon measures on compact spaces in the next section. 
Our aim is to emphasize that fundamental properties of such measures and the spaces $P(K)$ can be straightforwardly derived from a relatively simple result concerning  regular extensions of finitely additive set functions, all without delving into functional analysis. In particular, the Riesz representation theorem is not needed to show
 the compactness of $P(K)$; instead may be treated as its consequence. One has to admit that the proofs
 we give here are, at times,  sketchy but the general idea should be clear.

\section{Basic things}\label{bt}
The results recorded here originated in Marczewski's concept of a `compact measure' from \cite{Ma53} that was
developed by Pfanzagl and  Pierlo \cite{PP66}; see also 
Bogachev \cite[1.12]{Bo07}. The crucial Theorem \ref{re:2} is taken from
Bachman and Sultan \cite{BS80}.

\subsection{Regular extensions}\label{re}
We consider here an abstract space $X$, algebras of its subsets and finitely additive (nonnegative and finite)
measures. For any family $\cF$ of subsets of $X$, we write $\la \cF\ra$ for the algebra it generates.
Given an algebra $\cA$ and  a set function $\mu$ on $\cA$ we write, for any $Z\sub X$,
\[ \mu^\ast(Z)=\inf\{\mu(A): A\in\cA, A\supseteq Z\}\mbox{ and }
\mu_\ast(Z)=\sup \{\mu(A): A\in\cA, A\subseteq Z\}.
\]

\begin{lemma} \label{re:1}
If $\mu$ is finitely additive on an algebra $\cA$ then, for a given set  $Z\sub X$, the formula
\[ \ol{\mu}   ((A\cap Z)\cup (B\cap Z^c))=\mu^\ast(A\cap Z) +\mu_\ast(B\cap Z^c)  \]	
defines an extension of $\mu$ to a finitely additive measure on $\cA(Z)=\la \cA\cup\{Z\}\ra$
such that $\ol{\mu} (Z)=\mu^\ast(Z)$.
\end{lemma}	

\begin{proof}
First notice that, indeed,	every element of $\cA(Z)$ can be written as 
$(A\cap Z)\cup (B\cap Z^c)$ where $A,B\in\cA$. It is routine to check that
$\cA\ni A\mapsto \mu^\ast(A\cap Z) $ defines an additive set function on the trace of $\cA$ on $Z$; the same
holds for $\cA\ni A\mapsto \mu_\ast(A\cap Z) $. Thus $\ol{\mu}$ is finitely additive on $\cA(Z)$.

If $A\in\cA$ then $A=(A\cap Z)\cup (A\cap Z^c)$ so $\ol{\mu}(A)=\mu^\ast(A\cap Z)+\mu_\ast(A\cap Z^c)$
which is equal to $\mu(A)$, again  by routine calculations.

Finally, $Z=X\cap Z$ gives  $\ol{\mu} (Z)=\mu^\ast(Z)$.
\end{proof}

Saying that $\cL$ is a lattice, we mean that $\cL$ is a family of subsets of $X$  containing $\emptyset$, and  
closed under finite unions
and intersections. If $\cL$ is a lattice contained in an algebra $\cA$ then a finitely additive measure $\mu$ defined on $\cA$ is
said to be $\cL$-regular if
\[\mu(A)=\sup\{\mu(L): L\in\cL, L\sub A\}\mbox{ for every } A\in \cA.\]


\begin{theorem}(\cite[Theorem 2.1]{BS80}) \label{re:2}
	Let $\mathcal{A}$ be an algebra of subsets of $X$ and  let $\cL_0\subseteq \mathcal{A}$ be a lattice.
	Further assume that  $\mu$ is an $\cL_0$-regular, finitely additive measure on $\mathcal{A}$. 

Given any lattice $\cL\supseteq \cL_0$,  $\mu$ can be extended to an $\cL$-regular finitely additive measure on
$\la \cA\cup\cL\ra$.		
\end{theorem}

\begin{proof}
Consider first the case when $\cL$ is a lattice generated by $\cL_0$ and one additional set $Z$.
Note that  $\cL$ is then the family of sets of the form $(L\cap Z)\cup L'$ where $L,L'\in\cL_0$.

We consider $\ol{\mu}$ as in Lemma \ref{re:1}; it is enough to check that
	$\ol{\mu} $ is $\cL$-regular.  Indeed, for any $A\in\cA$ we have
	\[ \ol{\mu} (A\cap Z)=\mu^\ast (A\cap Z)=\sup\{\mu(L\cap Z): L\in\cL_0, L\sub A \}, \] 
	\[ \ol{\mu} (A\cap Z^c)=\mu_\ast (A\cap Z^c)=\sup\{\mu(L): L\in\cL_0, L\sub A\cap K^c \}, \] 
and those formulas yield the required regularity.

Now we only need to realize that  general case follows by transfinite induction: we can add new elements of
$\cL$ one by one.		
\end{proof}	

Given a  lattice $\cL$, we say that $\cL$  is countably compact if every countable subfamily of $\cL$	
with the finite intersection property has nonempty intersection.
The following seems to be pretty standard --- note that such  an extension theorem enables us to avoid the
classical Caratheodory method of constructing measures.

\begin{theorem}\label{re:3}
	Let $\mathcal{A}$ be an algebra of subsets of $X$,  $\cL\subseteq \mathcal{A}$ be a lattice
	and suppose that $\mu$ is an $\cL$-regular, finitely additive measure on $\mathcal{A}$. 

If $\cL$ is countably compact and closed under countable intersections then
$\mu$ can be extended to a countably additive $\cL$-regular measure on the $\sigma$-algebra $\sigma(\cA)$
generated by $\cA$.	
\end{theorem}

\begin{proof}
We first note the following.
\medskip

\noindent {\sc Claim 1.}
If $A_n\in\cA$ form a decreasing sequence  then for every $\eps>0$ the is a decreasing sequence
$L_n\in\cL$ such that $L_n\sub A_n$ and $\mu(A_n\sm L_n)< \eps$ for every $n$. 
\medskip

Indeed, for every $n$ choose $L_n'\in\cL$ such that $L_n' \sub A_n$ and $\mu(A_n\sm L_n')<\eps/2^n$, and define
$L_n=\bigcap_{k\le n} L_k'$. Then
\[ \mu(A_n\sm L_n)\le \sum_{k\le n} \mu(A_n\sm L_k')\le \sum_{k\le n} \mu(A_k\sm L_k')\le \sum_{k\le n} \eps/2^{k}<\eps,\]
as required.
\medskip

Claim 1 implies that, by countable compactness of $\cL$, 
$\mu$ is $\sigma$-smooth on $\cA$, that is $\lim_n\mu(A_n)=0$ whenever
$A_n\in\cA$ form a decreasing sequence with empty intersection. 
\medskip

\noindent {\sc Claim 2.} We have $\mu(\bigcap_n L_n )=\lim_n\mu(L_n)$ whenever $L_n\in\cL $
and $L_1\supseteq L_2\supseteq L_2\supseteq\ldots$
\medskip

Let $L=\bigcap_n L_n$; given $\eps>0$ find $Z\in\cL$ such that $Z\sub X\sm L$ and
$\mu(Z)>\mu(X\sm L)-\eps$. Then $\bigcap_n(L_n\cap Z)=\emptyset$ so, by countable compactness, 
 there is $k$ such that $L_k\cap Z=\emptyset$. Then
\[ \mu(L_k) \le \mu(X\sm Z)<\mu(L)+\eps, \]
which verifies the claim.
\medskip

Consider now the family $\Sigma$ of  sets $B\sub X$ for which,
for every $\eps>0$, there are $L, L'\in \cL$ such that 
\[ L\sub B, L'\sub X\sm B \mbox{  and } \mu(L)+\mu(L')>\mu(X)-\eps.\]
We put $\ol{\mu}(B)=\sup\{\mu(L): \cL\ni L\sub B\}$ for $B\in\Sigma$.
By the very definition, $\Sigma$ is closed under complements; it is easy to see that  
$\Sigma$ is closed under finite unions; consequently,   $\Sigma$ is an algebra of sets. 
Moreover, $\overline{\mu}$ is an $\cL$-regular,  additive set function on $\Sigma$ which is $\sigma$-smooth. 

Since $\Sigma\supseteq \cA$, it remains to check that
$\Sigma$ is a $\sigma$-algebra which amounts to verifying that
$\Sigma$ is also closed under decreasing intersections.

Take a decreasing sequence $B_n\in\Sigma $ and let $B=\bigcap_n B_n$. Writing $c=\lim_n \ol{\mu}(B_n)$, for
a given $\eps>0$ pick $k$ with $\ol{\mu}(B_k)<c+\eps$. It is clear that we can apply Claim 1 also to the sets $B_n$
to get a decreasing sequence of $L_n\in \cL$ such that 
\[ L_n \sub B_n \mbox{ and } \ol{\mu}(B_n\sm L_n)<\eps/2.\]
Write $L=\bigcap_n L_n$ and  find  $L'\in\cL$ satisfying  
\[ L'\sub X\sm B_k \mbox{  and } \mu(L')> \ol{\mu}(X\sm B_k)-\eps/2.\]
Using Claim 2 we conclude that
\[ \mu(L)+\mu(L')\ge \mu(L_k)+\mu(L')  -\eps; \]
since $L\sub B$ and $L'\cap B=\emptyset$, 
we have checked that $B\in\Sigma$ and  this finishes the proof.  
\end{proof}

\subsection{Baire and Borel measures}
As we have already declared, $K$ always stands for a compact Hausdorff space. By $C(K)$ we denote
the Banach space of real-valued continuous functions on $K$.
We,  of course,  define the Borel $\sigma$-algebra $Bor(K)$ as the one generated by all
open subsets of $K$. Let us first note the following easy consequences of inner-regularity of Radon measures.

\begin{lemma}\label{tau}
For any $K$ and $\mu\in P(K)$
\begin{enumerate}[(i)]
\item $\mu(B)=\inf\{\mu(U): U\supseteq B, U\mbox{\rm open}\}$ for every $B\in Bor(K)$;
\item if $V$ is the union of a family $\cU$ of open sets  which is closed under finite unions then
$\mu(V)=\sup\{\mu(U):U\in\cU\}$. 
\end{enumerate}
\end{lemma}

\begin{proof}
For $(i)$ apply inner regularity to $K\sm B$. Part $(ii)$ follows by inner-regularity and compactness.
\end{proof}

Typically, in a nonmetrizable space $K$ there is  another $\sigma$-algebra, often much  smaller than $Bor(K)$,  
with respect to which all the functions from $C(K)$ are measurable.
 A set $Z\sub K$ is called a zero set if $Z=g^{-1}[0]$ for some $g\in C(K)$. 
We shall write $\cZ_K$ for the family of all zero sets in $K$. The Baire $\sigma$-algebra $Ba(K)$ is defined
as the smallest one containing $\cZ_K$. Properties of zero sets
are discussed in Engelking \cite{En}; for properties of Baire sets (also in the general setting)  see Wheeler \cite{Wh83}.

We list below a few basic properties of zero sets; actually, they  hold in every normal topological space. 

\begin{lemma}\label{bab:1}
\begin{enumerate}
\item 	For every $g\in C(K)$ and a closed set $H\sub\er$, $g^{-1}[H] $ is a zero set.
\item $\cZ_K$ is a lattice closed under countable intersection.
\item For every $Z\in \cZ_K$ there are zero sets $Z_n$ such that $X\sm Z=\bigcup_n Z_n$.
\item For every closed $F\sub K$ and open $U\supseteq F$ there is $Z\in\cZ_K$ such that
$F\sub Z\sub U$.
\end{enumerate}		
\end{lemma}	

We can now discuss properties of countably additive Baire measures.
	
\begin{lemma}\label{bab:2}
Every finite measure on $Ba(K)$ is $\cZ_K$-regular.
\end{lemma}

\begin{proof}
	This can be checked as in the final step of the proof of \ref{re:3}. Namely,
	consider the family $\Sigma$ of those sets $B\sub K$ for which,
	for every $\eps>0$, there are $Z, Z'\in \cZ_K$ such that 
	\[ Z\sub B, Z'\sub K\sm B \mbox{ and } \mu(Z)+\mu(Z')>\mu(K)-\eps.\] 
	
	In a similar manner, we check that $\Sigma$ is a $\sigma$-algebra
	so it remains to note that $\cZ_K\sub \Sigma$ which is a consequence of Lemma \ref{bab:1}(3).
\end{proof}	

\begin{theorem}\label{bab:3}
Every finite measure $\mu$ on $Ba(K)$ can be uniquely extended to a Radon measure.		
\end{theorem}	

\begin{proof}
	The existence of an extension of $\mu$ to a Radon measure $\ol{\mu}$ on $K$ follows from
	Theorem \ref{re:2} and Theorem \ref{re:3} (where we let $\cL_0=\cZ_K$ and 
	 $\cL$ be the lattice of all closed sets).
	
	The uniqueness of such $\ol{\mu}$ follows from the observation that $\ol{\mu}$ must satisfy the formula
\[ \ol{\mu}(F)=\inf\{\mu(Z): Z\in\cZ_K, Z\supseteq F\}\]
for every closed set $F$ --- this is a consequence of Lemma \ref{bab:1}(4).	
\end{proof}	

A compact space $K$ is zero-dimensional if $\clop(K)$, the algebra of
closed-and-open  (briefly, clopen) sets in $K$ separates points of $K$.

\begin{corollary}\label{bab:5}
If $K$ is zero-dimensional then every finitely additive measure on $\clop(K)$ extends uniquely to a Radon measure.	
\end{corollary}

\begin{proof}
We again refer to Theorem \ref{re:2} (with $\cL_0=\clop(K)$) and apply Theorem \ref{re:3}.
The uniqueness follows from the fact that every closed subset of $K$ is an intersection of
clopen sets. 
\end{proof}

\begin{proposition}\label{bab:6}
	Given $\mu\in P(K)$, for every $B\in Bor(K)$ there is $A\in Ba(K)$ such that
$\mu(A\btu B)=0$.	
\end{proposition} 	

\begin{proof}
Using regularity  of $\mu$ we can assume that $B$ is $F_\sigma$; in turn,   it is enough
to approximate a given closed set $F$. By the formula used in the proof of \ref{bab:3},
there is $Z\in \cZ_K$ such that $Z\supseteq F$ and $\mu(Z)=\mu(F)$.
\end{proof}

\subsection{Product measures}\label{pm}
For any cardinal number $\kappa$ we can consider the Cantor cube  $2^\kappa$($=\{0,1\}^\kappa$)
equipped with the `standard' product measure
$\lambda_\kappa$ which is uniquely determined by the requirement that 
\[ \lambda_\kappa(\{x\in 2^\kappa: x(\alpha)=1\})=1/2 \mbox{ for every } \alpha<\kappa.\] 
What is the domain of $\lambda_\kappa$? 

For any set of coordinates $I\sub\kappa$ we denote the projection $2^\kappa\to 2^I$ by $\pi_I$.
Formally speaking, $\lambda_\kappa$ is defined on the product $\sigma$-algebra
consisting of all the sets of the form $A=\pi_I^{-1}[B]$ where $I\sub\kappa$ is countable
and $B\in Bor(2^I)$. For such a set we say that $A$ is determined by coordinates in $I$. 

The  product $\sigma$-algebra mentioned above is actually the Baire $\sigma$-algebra $Ba(2^\kappa)$.
To see why note first that, by very definition, every basic open set in $2^\kappa$ depends on a finite number of coordinates;
so does every clopen set. Now, for a continuous function $g:2^\kappa\to\er$ and closed $F\sub\er$,
$A=g^{-1}[F]$ is a $G_\delta$ subset of $2^\kappa$ so it can be written as $A=\bigcap_n C_n$ for some $C_n\in \clop(2^\kappa)$
and therefore  $A$ depends on countably many coordinates.

It is now clear that,  for $\kappa>\omega$,  $Ba(2^\kappa)$ is much smaller than $Bor(2^\kappa)$ --- no singleton in 
$2^\kappa$ is a Baire set.  This might indicate that to speak of   the `usual product measure' as a Borel measure
 we   should refer to  its unique Borel extension. 
However, here the situation is far simpler: $\lambda_\kappa$
is completion regular in the sense of the following theorem  which is a particular instance of Kakutani's theorem, see
Choksi and Fremlin \cite{CF79} for more information.

\begin{theorem}\label{kakutani}
For every $\kappa$ and for every $B\in Bor(2^\kappa)$ there are 
$B_1, B_2\in Ba(2^\kappa)$ such that $B_1\sub B \sub B_2$ and
$\lambda_\kappa(B_2\sm B_1)=0$.
\end{theorem}

\begin{proof}
Fix $\kappa$ and write $\lambda=\lambda_\kappa$ for simplicity.
Sets $B$ satisfying the assertion clearly form  a $\sigma$-algebra so it suffices to check that
every open set $U\sub 2^\kappa$ can be approximated from below and from above by Baire sets.
Finding $B_1$ is easy (actually in any compact space): take closed sets $F_n\sub U$ such that $\lambda(U\sm F_n)<1/n$;
for every $n$, by Lemma \ref{bab:1}, there is a zero set $Z_n$ such that $F_n\sub Z_n\sub U$. Then 
$B_1=\bigcup_n Z_n$ is the required Baire set. 

Now $B_1=\pi^{-1}_J[A]$ for some countable
$J\sub\kappa$ where $A$ is an $F_\sigma$ subset of $ 2^J$.
Consider
\[ V=\pi_{J}^{-1}[U'] \mbox{ where } U'=\pi_J [U];\] 
here $U'$ is open in $2^J$ so $V$ is an open Baire set containing $U$ and it remains to check that it may play a role of $B_2$, that is
$\lambda(V)=\lambda(U)$. 

 Write $\lambda'=\pi_J[\lambda]$ and
$\lambda''=\pi_{\kappa\sm J}[\lambda]$
for the image measures --- see \ref{bppk} for the notation. Then $\lambda=\lambda'\otimes\lambda''$, i.e.\
we may  think that $\lambda$ is a product measure on $2^\kappa=2^J\times 2^{\kappa\sm J}$.
Suppose that $\lambda (V\sm U)>0$. Then $\lambda'(U'\sm A)>0$ so there is $x'\in U'\sm A$ with the property that
$\lambda'(W'\cap (U'\sm A))>0$ for every open $W'\sub 2^J$ neighbourhood of $x'$.
As $x'\in U'=\pi_J[U]$, there is $x''\in 2^{\kappa\sm J}$ such that $x=(x',x'')\in U$.
But $U$ is open so $x\in W'\times W''\sub U$ for some basic open sets $W',W''$ We arrive at contradiction
examining the set $Z=(W'\sm A)\times W''$: On one hand,
\[\lambda (Z)=\lambda'(W'\sm A)\lambda''(W'')>0.\] 
On the other hand, $\lambda(Z)=0~$ because $Z\sub U\sm B_1$.  
\end{proof}

 Let us recall that  the product $\sigma$-algebra coincides with the Baire $\sigma$-algebra
 in an arbitrary product of separable metrizable spaces, see e.g.\ \cite{CN82}.
 The argument of Theorem \ref{kakutani} can be easily adapted to the product spaces $[0,1]^\kappa$
 with the product of the Lebesgue measure on $[0,1]$.

\subsection{Spaces of measures}\label{bppk}
First of all we note that if $K$ is zero-dimensional then, by Corollary \ref{bab:5},
$P(K)$ may be identified with $P(\clop(K))$, the space of all finitely additive probability measures
on the algebra of sets $\clop(K)$. 
This identification observes the topological structure on $P(K)$ when we simply equip 
$P(\clop(K))$ with the topology  inherited from $[0,1]^{\clop(K)}$, that is
the topology of convergence on clopen sets. Indeed, convergence on $\clop(K)$ 
yields  convergence on simple continuous functions and
these form a  norm-dense subspace of  $C(K)$.  As $P(\clop(K))$ is a closed subset of  $[0,1]^{\clop(K)}$,
we conclude the following.

\begin{lemma}\label{bp:1}
The space $P(K)$	is compact for every zero-dimensional $K$.	
\end{lemma}	

If $f:K\to L$ is a continuous surjection between compact spaces then for every $\mu\in P(K)$
we write $f[\mu]$ for the image measures defined on $L$, that is
\[ f[\mu](B)=\mu(f^{-1}[B])\mbox{ for every } B\in Bor(L). \]
Note that $\nu=f[\mu]\in P(L)$, i.e.\ $\nu$ is regular:
for any $B\in Bor(L)$ we have $B'=f^{-1}[B]\in Bor (K)$ so, for a given $\eps>0$,  there is
a closed set $F\sub B'$ such that  $\mu(B'\sm F)<\eps$. Then $H=f[F]$ is a 
closed subset
of $B$ and $F\sub f^{-1}[H]\sub B'$ so $\nu(B\sm H)<\eps$.

\begin{theorem}\label{bp:2}
If $f:K\to L$ is a continuous surjection then
\[ P(K)\ni\mu\mapsto f[\mu]\in P(L)\]
is a continuous surjection between the corresponding spaces of measures.	
\end{theorem}

\begin{proof}
To check surjectivity fix $\nu\in P(L)$; to find $\mu\in P(K)$ such that $f[\mu]=\nu$ consider
\[ \cA=\{f^{-1}[B]: B\in Bor(L)\},\quad \cL_0=\{f^{-1}(H):H=\overline{H}\subseteq L \}.\]
Then $\cA$ is a $\sigma$-algebra of sets, $\cL\sub \cA$ is a lattice and putting
$\mu_0(f^{-1}[B])=\nu(B)$ for $B\in Bor(L)$ we have an $\cL_0$-regular measure $\mu_0$ defined on $\cA$.
By Theorems \ref{re:2} and \ref{re:3}, $\mu_0$ extends to the required $\mu\in P(K)$.

The continuity of the mapping $\mu\mapsto f[\mu]$ follows by  the standard formula
\[ \int_L g\; {\rm d}f[\mu]=\int_K g\circ f\; {\rm d}\mu\]
satisfied by any $g\in C(L)$.
\end{proof}

Let us mention here that Theorem \ref{bp:2} implies that every compact space $K$ that can be continuously
mapped onto $[0,1]$ (i.e.\ every $K$ which is not scattered) carries a nonatomic probability measure
(here: the one vanishing on points).

\begin{theorem}\label{bp:3}
The space $P(K)$ is compact for every $K$.
\end{theorem}

\begin{proof}
Recall that every compact space is a continuous image of a compact zero-dimensional space;
then apply Lemma \ref{bp:1}   and Theorem \ref{bp:2}.
\end{proof}

\begin{remark}\label{add:1}
In the setting of Theorem \ref{bp:2}, for a given $\nu\in P(K)$ we can find
$\mu\in P(K)$ such that $f[\mu]=\nu$ and, moreover,
the $\sigma$-algebra $\Sigma=\{f^{-1}[B]: B\in Bor(L)\}$ is $\bigtriangleup$-dense in
$Bor(K)$, that is for every $A\in Bor(K)$ there is $E\in\Sigma$ such that
$\mu(A\bigtriangleup B)=0$.

To see this, one has to examine briefly  the proof of Theorem \ref{re:2}. Alternatively, we can
recall a general result due to Douglas \cite{Do64}. Namely,
the set $M=\{\mu\in P(K): f[\mu]=\nu\}$ is a compact convex subset of $P(K)$
and any extreme point of $M$ is a measure having the required property.
\end{remark}

Most often, Theorem \ref{bp:3} is derived from the Banach-Alaoglu theorem on $weak^\ast$ compactness
in dual Banach spaces; this
however, requires the Riesz representation theorem. Note that we could  actually do the opposite, derive the Riesz theorem from
\ref{bp:3} using the observations given below.
 
 \newcommand{\co}{\protect{\rm conv}}
 
 For every $K$ the space $P(K)$ contains the subspace of Dirac measures 
 $\Delta_K=\{\delta_x: x\in K\}$. Its convex hull $\co(\Delta_K)$ is the subset of $P(K)$ consisting of
 all finitely supported probabilities. 
 
  \begin{lemma}\label{bp:4}
  	The set $\co(\Delta_K)$ is dense in $P(K)$ for every $K$.
  \end{lemma}	
  
 \begin{proof}
 Rather	 than proving this directly, note first that the result is almost immediate for zero-dimensional $K$:
 for any $\mu\in P(K)$ and a finite partition $K=\bigcup_ i A_i$ into nonempty clopen set, 
 choose any $x_i\in A_i$ and observe that
 the measure $\nu=\sum_i \mu(A_i)\delta_{x_i}$ agrees with $\mu$ on $A_i$'s.
 
 For the general case use the fact that $P(K)$ is a continuous image of the space of probabilities
 on a zero-dimensional compactum.	
 \end{proof}	 
 
 If $\mu\in P(K)$ then the support of $\mu$ is the smallest closed set $S\sub K$ with the property
 $\mu(S)=1$. The existence of such $S$ follows from the observation that
  $K\sm S$ is the union of all open sets of measure zero and $\mu(K\sm S)=0$ by Lemma \ref{tau}(ii).
  When  $S=K$ then we say that $\mu$ is strictly positive which amounts to stating that
  $\mu(U)>0$ for every nonempty open set $U\sub K$. 
   Note  also that for a closed set $S\sub K$, the space $P(S)$ may be treated as a subspace of $P(K)$ consisting of measures
   concentrated on $S$. 
 
 \subsection{Kelley's intersection number}
 Kelley introduced in \cite{Ke59}   that concept to characterize  Boolean
 algebras carrying strictly positive measures. For a family $\cF$  of sets, the intersection number
 of $\cF$, denoted here by $i(\cF)$,  can be defined by saying that 
 $i(\cF)\ge c$ if for every $n$ and any $F_1,\ldots, F_n\in\cF$ we have
 \[ \frac{\|\sum_{k\le n}\chi_{F_k}\|}{n}\ge c.\]
 Here $\|\cdot\|$ stands for the supremum norm  so the fraction above
 indicates the proportion of sets in a given sequence having nonempty intersection.
 
 Intersection numbers enables one to express the property that a given compact space $K$ carries a strictly
 positive measure $\mu$  as a kind of chain condition, see
 Comfort and Negrepontis \cite[Chapter 6]{CN82} and \stevo\ \cite{To00}. 
 
 We shall need the following result which  is a particular case of \cite[Lemma 6.3]{CN82}; see also
 M\"agerl and Namioka \cite[Theorem 1]{MN80} or Kindler \cite{Ki88}.
 We sketch another proof which indicates that the essence of Kelley's theorem can be reduced,
in  the language of   linear programming, to duality. 
 
 \begin{theorem} \label{kelley} 
 If $\cF$ is a family of closed subsets of $K$ and  $i(\cF)\ge c$  then 
 there is $\mu\in P(K)$ such that $\mu(F)\ge c$ for every $F\in\cF$. 
 \end{theorem}
 
 \begin{proof}
Note first that for any family
$\cH$ of closed subsets of $K$, the set
\[ P(\cH,c)=\{\mu\in P(K): \mu(H)\ge c\mbox{ for every } H\in\cH\}\]
is closed in $P(K)$. Moreover, $P(\cF,c)$ is the intersection of $P(\cH,c)$
over all finite $\cH\sub\cF$. Hence, by compactness of $P(K)$,  is suffices to prove the assertion for
a finite family $\cF$. 

Consider such a finite family  $\cF=\{F_1,\ldots, F_k\}$.
Then the algebra $\cA$ generated by $\cF$ is also finite; let $T_1,\ldots, T_n$ be the list of all atoms
of $\cA$. Pick $p_i\in T_i$ for every $i\le n$. We shall prove that there is $\mu\in P(K)$ of the form
$\mu=\sum_{i\le n} x_i\delta_{p_i}$ such that $\mu(F_j) \ge c$ for $j\le k$.

Writing $A_j=\{i\le n: p_i\in F_j\}$ we see that the problem reduces to
 finding nonnegative reals $x_i$ such that 
 \[ \sum_i x_i=1 \mbox{ and } \sum_{i\in A_j} x_i\ge c \mbox{  for every } j\le k.\]
  This may be seen as an optimization problem:
 find $\min \sum_i x_i$ subject to the above constraints.

Recall here that a typical  problem in linear programming is to find 
\[ \min \sum_{i\le n} a_i x_i \quad \mbox{  subject to }\quad  Mx\ge b, x\ge 0,\]
where  $x=(x_1,\ldots, x_n)\in \er^n$,  $M$ is an $k\times n$ matrix while $b\in \er^k$.
 Then the dual problem is: 
 \[ \max \sum_{j\le k} b_jy_j \quad \mbox{ subject to }\quad  A^{T}y\le a, y\in \er^k, y\ge 0.\]
The main use of the duality is that if $x$ and $y$ are {\em feasible} for the the original problem and the dual one, respectively, then $\sum_{j\le k} b_jy_j\le \sum_{i\le n} a_ix_i$; moreover, the equality holds here
if $x$ and $y$ are optimal; see  Bertsimas and  Tsitsiklis \cite[Theorem  4.4]{BT97}.

Imagine now a matrix $M$ coding our requirements  $\sum_{i\in A_j} x_i\ge c$
and observe that the constraints $M^Ty\le 1$ in the dual problem mean that
$\|\sum_{j\le k} y_j\chi_{A_j}\|\le 1$. Then note that $i(\{A_1,\ldots, A_k\})\ge c$ implies that
\[  \|\sum_{j\le k} y_j\chi_{A_j}\|\ge c\sum_{j\le k}y_j\]
 for any nonnegative $y_j$'s.
We conclude that the values of the objective function in the dual problem is bounded
from above by 1 hence the original problem allows a solution such that $\sum_ix_i\le 1$, and
we are done.
\end{proof}

 \section{Boolean algebras and some set theory}
 Given a Boolean algebra $\fB$,  $\ult(\fB)$ denotes its Stone space, i.e.\  the set of all ultrafilters on $\fB$
 endowed with the topology having as a base  the collection of sets of the form
 $\widehat{B} = \{p\in \ult(\fB)\colon B\in p\}$, for  $B\in \fB$. Then $\ult(\fB)$ is a compact Hausdorff space and the assignment $B\to \widehat{B}$ defines a Boolean
 isomorphism between $\fB$ and the algebra of clopen subsets of $\ult(\fB)$.
 Conversely, every compact Hausdorff space possessing a base of clopen sets (i.e. a $0$-dimensional space) 
 can be seen as the Stone space of its algebra of clopen sets; 
 Koppelberg \cite{kopp} is our basic reference here. 
 
 Most often, we consider algebras of sets rather than abstract Boolean algebras.
 Following the notation of subsection \ref{bppk}, if $\cA$ is an algebra of sets then 
 $P(\cA)$ stands for the space of all {\em finitely} additive probability measures on $\cA$. 
 Recall that $P(\cA)$ is always equipped with the topology inherited from $[0,1]^\cA$.
  
 \begin{lemma}\label{xy}
 For every algebra $\cA$ the space $P(\cA)$ is homeomorphic to $P(K)$ where
 $K=\ult(\cA)$.
 \end{lemma}
 
 \begin{proof}
 If $\mu\in P(\cA)$ then putting 
 $\wh{\mu}(\wh{A})=\mu(A)$ for $A\in \cA$
 we define a finitely additive measure $\wh{\mu}\in P(\clop(K))$ which, by Corollary \ref{bab:5},
 can be uniquely extended to an element of $P(K)$.
  \end{proof}

 \subsection{Measure algebras and Maharam types}
 For every finite measure space $(X,\Sigma,\mu)$, writing  
 $\cN=\{A\in\Sigma: \mu(A)=0\}$,
 one can form its measure algebra
 $\fA=\Sigma/\cN$ of all equivalence classes $A\kr=\{B\in\Sigma: \mu(B\btu A)=0\}$, $A\in \Sigma$.
 Moreover, we can treat $\mu$ as being defined on $\fA$.
 Measure algebras are discussed in Fremlin \cite{Fr89} and \cite{Fr4}.

 For every $\kappa$, write $\fB_\kappa$ for the measure algebra of the measure $\lambda_\kappa$ on $2^\kappa$.
  By  completion  regularity of product measures (see subsection \ref{pm}),  there is no difference
  if we form $\fB_\kappa$ for the measure $\lambda_\kappa$ defined on
  $Ba(2^\kappa)$ or with respect to its unique regular  extension to $Bor(2^\kappa)$. 
 This explains why $\fB_\kappa$ is of size $\kappa^\omega$ (for $\kappa>1$).

 \begin{definition}\label{mt}
 Given $\mu\in P(K)$, the Maharam type $\mt(\mu)$ of  $\mu$ is  the minimal infinite cardinal number
 $\kappa$ such that  there is  $\cD\sub Bor(K)$ with $|\cD|=\kappa$ and having the property  that 
 for any $\varepsilon>0$ and any $B\in Bor(K)$ there is $D\in\cD$ such that
  $\mu(D\bigtriangleup B)<\varepsilon$. 
  \end{definition}
  
  In other words, $\mt(\mu)$ is the density of the corresponding measure algebra $\fA$
  with respect to the metric $\fA\times\fA\ni (a,b)\to \mu(a\btu b)$. It is also not difficult to check that 
  $\mt(\mu)$ is equal to the density of the Banach space $L_1(\mu)$.
  
  Recall that, by the Maharam structure theorem, if $\mu$ is of type $\le\kappa$ then
  its measure algebra $\fA$ can be embedded into $\fB_\kappa$ by `the
  measure preserving Boolean homomorphism'.
  Moreover, if $\mu$ is homogeneous of type $\kappa$ (i.e.\ its type
  is $\kappa$ when we restrict $\mu$ to any set of positive measure) then
  $\fA$ is isomorphic to $\fB_\kappa$, see e.g.\ Fremlin \cite{Fr89}.   
  
  \begin{remark}\label{add:2}
  It follows immediately from `Remark \ref{add:1} that if $f:K\to L$ is a continuous surjection, then
  for a given $\nu\in P(L)$,  there is $\mu\in P(K)$ such that $f[\mu]=\nu$ and
  $\mu$ has the same Maharam type as $\nu$.
  \end{remark}

  It will be convenient to use the following notion.
  
  \begin{definition}
  We say that a compact space $K$ is in the class \MS\ if $K$ carries only measures of type $\omega$. 
  \end{definition}
  
  This terminology  is borrowed from D\v{z}amonja and Kunen \cite{DK95}; MS stands for {\em Measure Separable} --- measures of countable type are also simply called  separable.
  The class \MS\ contains all metric compacta and is closed under closed subspaces, continuous images and countable products, see
  \cite{DK95} for more information.

   \subsection{Calibers of measures}\label{cm}
 
 We shall use the following concept.
 
 \begin{definition}\label{cm:1}
 A  cardinal number  $\kappa$ is a {\em caliber}  of a measure $\mu\in P(K)$  if, for every
 family $\{A_{\alpha}:\alpha<\kappa\}$ in
 $Bor(K)$ of $\mu$-positive sets, there is $x\in K$ such that
 $\{\alpha<\kappa: x\in A_\alpha\}$ has cardinality $\kappa$.
 \end{definition}

 Basic properties of calibers of measures are discussed in \cite{DP04} and Fremlin \cite[525]{Fr5}.
 In particular, \cite[Lemma 2.5]{DP04} states that a given cardinal number $\kappa$ is
 a caliber of all Radon measures if $\kappa$ is a precaliber of all measure algebras, that
 is if $\{a_\xi:\xi<\kappa\}$ is a family of positive elements of such an algebra then there
 is as set $I\sub\kappa$ of size $\kappa$ such that $\{a_\xi:\xi\in I\}$ is centered. 
 By the Maharam theorem, this yields the conclusion that $\kappa$ is a caliber of Radon measures if and only if
 $\kappa$ is a caliber of the measure $\lambda_\kappa$ on $2^\kappa$.
 
 For the following basic fact see \stevo\ \cite[Lemma 6]{To96} (or \cite[525G]{Fr5}, \cite[Section 3]{DP04}).
 
\begin{lemma}\label{cm:2}
The cardinal number $\omega_1$ is a caliber of Radon measures if and only if
 the Cantor cube $2^{\omega_1}$ cannot be covered by $\omega_1$ many sets
that are negligible with respect to $\lambda_{\omega_1}$.
 \end{lemma}

In the language of cardinal coefficient of ideals, see \cite[Chapter 54]{Fr5}, the condition mentioned  above can be written as 
  ${\rm cov} (2^{\omega_1})>\omega_1$. Recall that 
 ${\rm cov} (2^{\omega_1})\le  {\rm cov} (2^{\omega})$ and the sharp inequality is relatively consistent, see 
 Kraszewski \cite[Remark after Theorem 5.5]{Kr01}.

   Concerning Martin's axiom, the reader is referred to Jech \cite[Chapter 16]{Jech}
  and Fremlin \cite{Fr84}. 
Here we write $\MA(\kappa)$ for Martin's axiom stated for
 $\kappa$ many dense sets in a $ccc$ poset and $\MA$
 means that $\MA(\kappa)$ holds for every $\kappa<\con$.  
 
 Note that, by Lemma \ref{cm:2}, $\omega_1$ is not a caliber of Radon measures if
 the continuum hypothesis $\CH$ holds. On the other hand,
 $\MA(\omega_1)$ implies that $\omega_1$ is a caliber of measures.

\section{Cardinality of $P(K)$ and separability }\label{btp}

For every $K$ the space $P(K)$ contains $\Delta_K=\{\delta_x: x\in K\}$ which is clearly homeomorphic to $K$.
It follows that if $P(K)$ has some topological property which is hereditary with respect to closed subspaces then
the space $K$ enjoys the same property.
It is easy to check that $P(K)$ is metrizable if and only if $K$ is metrizable (recall that metrizability of $K$ is
equivalent to separability of the Banach space  $C(K)$). As we shall see below, there is no such straightforward relation
between $K$ and $P(K)$ when we inspect other  topological properties.

\subsection{Cardinality}\label{card}
 For every $K$ we have
 \[ |K|^\omega\le |P(K)|\le 2^{w(K)};\]
here $w(\cdot)$ denotes the topological weight. Indeed, 
$|K|^\omega$ is the number of fully atomic measures on $K$. For the upper bound above note that
 if we take a base $\cU$ for the topology of $K$ with $|\cU|=w(K)$ then we can assume that $\cU$ is closed under finite unions. 
 Now for any $\mu\in P(K)$ and any open set $V$ we have 
 \[ \mu(V)=\sup\{\mu(U): U\in\cU, U\sub V\},\]
  by Lemma \ref{tau}(ii). In turn,  by Lemma \ref{tau}(i), 
  $\mu$ is determined by its values on $\cU$. 
  
  Lipecki \cite{Li01} asked whether
  actually $|P(K)|=|K|^\omega$ for every compact space $K$. 
  
  \begin{theorem}[Fremlin and Plebanek \cite{FP03}]
  Under \MA\ there is a compact space $K$ of cardinality $\con$ that carries $2^\con$ many pairwise
  orthogonal Radon measures (in particular,  $|P(K)|=2^\con$).
  \end{theorem}
 
 However, the following problem seems to be open.
 
 \begin{problem}
 Is there a ZFC example of a compact space of cardinality $\con$ for which $|P(K)|>\con$?
 \end{problem}
 
 Dales and Plebanek \cite{DP19} examine  questions related to cardinality of $P(K)$ in connection with
 the structure of dual Banach spaces of the form $C(K)^\ast$.
 
 \subsection{Separability}\label{sep}
 Arguing as for the proof of Lemma \ref{bp:4} we can check that  if $D$ is a dense subset of $K$ then $P(K)$ is the closed  convex hull of the
 set $\{\delta_x:x\in D\}$. It follows that $P(K)$ is separable for every separable compactum $K$. As we shall see,
 this cannot be reversed and the question {\em when $P(K)$ is separable?} is delicate.

\begin{lemma}\label{approx}
The space $P(K)$ is separable if and only if there is a sequence of $\mu_n\in P(K)$ such that
for every nonempty open set $U\sub K$ there is $n$ such that $\mu_n(U)>1/2$. 
\end{lemma}

The above observation was used by Talagrand \cite{Ta80} and
M\"{a}gerl and Namioka \cite{MN80} who showed that the separability of $P(K)$ is
equivalent to some chain-type condition for the topology of $K$ (here $1/2$ can be replaced by any
constant $c\in (0,1)$).

\begin{corollary}[\cite{Ta80}]\label{corta}
If a compact space $K$ carries a strictly positive measure $\lambda$ of countable type then
$P(K)$ is separable.
\end{corollary}  

\begin{proof}
Take a countable family $\cA$ in $Bor(K)$ such that 
$\inf\{\lambda(B\btu A): A\in \cA\}=0$
for every Borel set $B$. Then let $\mu_A$ be the normalized restriction of $\lambda$
to $A$ (whenever $\lambda(A)>0$). Then apply Lemma \ref{approx} to
the family of $\mu_A$, $A\in\cA$.
\end{proof}
 
 Now it is clear that if we consider the Stone space $K$ of the Lebesgue measure algebra
 then $P(K)$ is separable while $K$ is not.
 Lemma \ref{corta} gives a handy sufficient condition for separability of $P(K)$ which is, however,
 not necessary. The following was proved by Talagrand \cite{Ta80} under CH, and
 by D\v{z}amonja and Plebanek \cite{DP08} in ZFC. 
 
 \begin{theorem}
 There is a compact space $K$  such that $P(K)$ is separable but $K$
 carries no strictly positive measure of countable type.
 \end{theorem}
 
 It is worth noting here that the separability of $P(K)$ 
 is equivalent to  $weak^\ast$ separability of the unit ball $M_1(K)$ in $C(K)^\ast$. In turn, the latter property
 implies that $M(K)=C(K)^\ast$ is $weak^\ast$ separable but not vice versa as it was proved, again, 
 by Talagrand \cite{Ta80} under CH and later by
 Avil\'es, Plebanek and Rodr\'{\i}guez  \cite{APR14} in ZFC.

 \section{First countability and Corson compacta}\label{fccc}
 A compact space $K$ is {\em Corson compact} if $K$
 can be embedded into a $\Sigma$-product $\Sigma([0,1]^\Gamma)$ for some index set $\Gamma$;
 recall that  
 \[ \Sigma([0,1]^\Gamma)=\{x\in [0,1]^\Gamma: |\{\gamma: x(\gamma)\neq 0\}|\le\omega\}.\] 
 Corson compacta
 form an important class of spaces related to functional analysis, see
 Argyros, Mercourakis and Negrepontis \cite{AMN}
 and Kalenda \cite{Ka00}.
 Corson compacta appear here in the company of first-countable spaces since
 the measures on spaces from both the classes behave, to some extend,  in a similar manner.

 \begin{definition}
A measure $\mu\in P (K )$ is said to be {\em strongly countably determined}
 if there exists a continuous map $f: K\to [0, 1]^\omega$ such
that for any compact set $F\sub K$ we have 
$\mu(F)=\mu(f^{-1}[f[F]])$.
\end{definition}

Strongly countably determined measures were introduced
by Babiker [1] who called them  {\em uniformly regular measures}.

\begin{lemma} [\cite{Ba76}] \label{sct}
A measure $\mu\in P(K)$ is strongly countably determined if and only if 
there is a countable family $\cZ$ of zero subsets of K such
that
\[\mu(U)=\sup\{ \mu(Z): Z\in\cZ, Z\sub U\},\]
for every open set $U\sub K$.
\end{lemma}

The latter condition seems to justify the name strongly countably determined; note that
if we relax the property to saying that $\cZ$ is a countable family of closed sets then we get a weaker property and
the measure $\mu$  in question is called  then
{\em countably determined}. The following implications are obvious 
\[ \mbox{strongly countably determined } \Rightarrow
\mbox{ countably determined } \Rightarrow \mbox{ of countable type}.\]
The next result  is due to Pol \cite[Proposition 2]{Po82}.

\begin{proposition}\label{fact:Pol}
For any compact space $K$ and a measure $\mu\in P(K)$,
the space $P(K)$ is first-countable at $\mu$ if and only if $\mu$ is strongly countably determined.
\end{proposition}

The proposition above can be generalized to its higher cardinal analogue, see
Krupski and Plebanek \cite{KP11} also for the proof of the next result.

\begin{theorem}\label{krupski}
Given a compact space $K$ in the class \MS, the space
$P(K)$ has a dense subset of $G_\delta$-points.
\end{theorem} 

The next theorem summarizes the results due to Kunen and van Mill \cite{KM95}
Argyros, Mercourakis and Negrepontis \cite{AMN}
and Plebanek  \cite{Pl95a}.

\begin{theorem}\label{all}
The following are equivalent
\begin{enumerate}[(i)]
\item $\omega_1$ is a caliber of Radon measures;
\item whenever $K$ is Corson compact and $\mu\in P(K)$  then $\mu(K_0)=1$ for some
 metrizable closed set $K_0\sub K$;
 \item every Corson compact is in \MS;
  \item $P(K)$ is Corson compact for every Corson compact $K$;
   \item every first-countable compactum $K$ is in the class \MS.
\end{enumerate}
\end{theorem}

\begin{proof}
$(i)\to (ii)$. Suppose that $K\sub\Sigma([0,1]^\Gamma)$, take $\mu\in P(K)$  and
write
\[ \Gamma_0=\{\gamma\in\Gamma: \mu( V_\gamma)>0\}, \mbox{ where } 
V_\gamma=\{x\in K: x(\gamma)>0\}. \]
If we suppose that $\mu$ does not satisfy $(ii)$ then
we conclude that   the set $\Gamma_0$ is uncountable.  
But no family of uncountably many $V_\gamma$ has nonempty intersection by
the very definition of $\Sigma$-product, and it follows that $\omega_1$ is not a caliber of the measure $\mu$.

$(ii)\to (iii)$ follows from the fact that
every measure concentrated on a metrizable compactum has type $\omega$.

To prove $(iii)\to (i)$ suppose that $\omega_1$ is not a caliber of Radon measures.
By Lemma \ref{cm:2},   $\omega_1$ is not a caliber of the measure $\lambda_{\omega_1}$ on the Cantor cube $2^{\omega_1}$.
Say that this is witnessed by a family $\{B_\xi:\xi<\oo\}$ of sets of positive measure.

We inductively choose closed sets $F_\xi\sub 2^\oo$ of positive measure so that for every $\xi<\oo$

\begin{enumerate}[(i)]
\item $F_\xi\sub B_\xi$;
\item $c_\xi=\inf\{ \lambda(F_\xi\btu A): A\in \cA_\xi\}>0$, where $\cA_\xi$ is the algebra of sets
generated by $\{F_\beta: \beta<\xi\}$. 
\end{enumerate} 
 
 The construction is straightforward as $\lambda_\oo$ is regular and has type $\oo$ on every set of positive measure.
  
 Then we consider the algebra $\cA$ of subsets of $2^\oo$ generated by the whole family
 $\{F_\beta: \beta<\oo\}$ and check that $K=\ult(\cA)$ is a Corson compact space carrying a measure
 of type $\oo$. Write $\mu$ for $\lambda_\oo$ restricted to $\cA$. Then $\mu\in P(\cA)$ so,
 by Lemma \ref{xy}, $\mu$ can be seen as a Radon measure on $K$. Condition $(ii)$ above assures that
 $\mu$ is of uncountable type: indeed,  there is $c>0$ such that $c_\xi\ge c$ for uncountable many $\xi$. 
 We get an uncountable set $I\sub\oo$ such that $\mu(F_\xi\btu F_\eta)\ge c$ whenever $\xi\neq\eta$ are in $I$
 which means that $\mu$ is of type $\oo$.  
 
 To see that $K$ is Corson compact consider the diagonal mapping 
  \[ g=\mathop{\BTU}_{\xi<\oo} \chi_{\wh{F_\xi}}: K\to 2^\oo;\]
 $g$ us one-to-one (if two ultrafilters in $\ult(\cA)$ agree on all the generators $F_\xi$ then they are identical).
 It follows that $g$ is an embedding and in takes values in $\Sigma(2^\oo)$ by $(i)$.
 
 $(iii)\leftrightarrow (iv)$ can be found in \cite[Section 3]{AMN}; recall that a Corson compact space is metrizable if and only if it is separable.
 
 Finally, $(i)\to (v)$ is the main result of \cite{Pl95a} while
 $(v)\to (i)$ was first proved in \cite{KM95}; another construction, in the spirit of that above,
  can be found in \cite[Theorem 5.1]{Pl97}.
\end{proof}

Let us  mention here that Theorem \ref{all} has higher-cardinals analogues, see
\cite[section 12]{MPZ} where the so called $\kappa$-Corson compact spaces are
discussed.

Core ideas of some \CH\ constructions related to Theorem \ref{all} were invented by  
Haydon \cite{Ha78}, Talagrand \cite{Ta80} and Kunen \cite{Ku81}.
Their objective were different but their constructions had  common ingredients
which were later developed in various directions, see e.g.\ 
D\v{z}amonja and Kunen \cite{DK93}, Plebanek \cite{Pl23},
Banakh and Gabryieylan \cite{BG24}, Koszmider and Silber \cite{KS24}.

When we consider first-countable compact spaces or Corson compacta $K$ and
examine the corresponding $P(K)$ space then \CH\ and $\MA(\oo)$ yield completely  different pictures.
 We shall illustrate this phenomenon proving the following
result due to Talagrand \cite{Ta80a}.

\begin{theorem}\label{underch}
Under \CH\ there is a first-countable Corson compact space $K$ such that
$P(K)$ contains a topological copy of $\beta\omega$.
\end{theorem}

\begin{proof}
Using \CH\ we fix an enumeration $\{t_\xi:\xi<\oo\}$ of the Cantor set $2^\omega$
and list all the subset of $\omega$ as $\{N_\xi:\xi<\oo\}$. For every $\xi$ we choose
a closed set $F_\xi\sub 2^\omega\sm \{t_\beta: \beta<\xi\}$ such that $\lambda(F_\xi)>1/2$
(here $\lambda=\lambda_\omega$ is the usual product measure on the Cantor set). 
Consider two families  of subsets of $2^\omega\times\omega$:
\[\cF=\{F_\xi\times N_\xi:\xi<\oo\} \mbox{ and } \cC=\{C\times\omega: C\in\clop(2^\omega)\}; \]
we check below that our target space $K$ will be $K=\ult(\cA)$ where $\cA$ is
the algebra of subsets of $2^\omega\times\omega$ generated by $\cF\cup\cC$.

Take any ultrafilter $\pe\in\ult(\cA)$; then there is a unique $t\in 2^\omega$ such that
for any clopen $C\sub 2^\omega$ we have
\[ C\times\omega \in p \mbox{ if and only if } t\in C.\]
Then $t=t_\alpha$ for some $\alpha<\oo$ and it is enough to check that the countabe algebra 
\[ \cA_0=\la \{F_\beta\times\omega:\beta\le \alpha\}\cup\cC\]
  contains a base of the ultrafilter $p$, 
that is, for every $A\in p$ there is $A_0\in\cA_0$ such that $A_0\in p$ and $A_0\sub A$ (which immediately translates to
the fact that $\ult(\cA)$ is first-countable at $p$). This follows from the following observation.

For any $\xi>\alpha$ we have $t\notin F_\xi$ so there is a clopen set $C\sub 2^\omega$
such that $t\in C$ and $C\cap F_\xi=\emptyset$. Then $C\times\omega\in p$
and $C\times\omega\sub 2^\omega\times\omega\sm (F_\xi\times B_\xi)$.   

The fact that $K=\ult(\cA)$ is Corson compact can be checked as in the proof of  Theorem \ref{all};
simply note that  every $p\in\ult(\cA)$ can contain only a countable number of generators.

For a given $n\in \omega$, we can define a finitely additive measure
$\mu_n\in P(\cA)$ by the formula $\mu_n(A)=\lambda(A^n)$ for $A\in\cA$; here
$A^n=\{t\in 2^\omega: (t,n)\in A\}$ denotes the horizontal section of $A$. 
Using Lemma \ref{xy} we may think that $\mu_n\in P(K)$ and it suffices to
check that the closure of the set $\{\mu_n: n\in\omega\}$ inside $P(K)$
is homeomorphic to $\beta\omega$. This amounts to showing that
\[ (*)\quad \ol{ \{\mu_n: n\in N\}}\cap \ol{ \{\mu_n: n\in \omega\sm N\}}=\emptyset\]
for every $N\sub \omega$.

Given $N\sub\omega$, we have $N=N_\xi$ for some $\xi<\oo$. Then
\[ \mu_n(F_\xi\times N_\xi)=\lambda(F_\xi)>1/2 \mbox{ whenever }n\in N,\] 
while the value above is  $0$ for $n\in\omega\sm N$.
This mean that the clopen set $\wh{F_\xi\times N_\xi}$ indicates that $(*)$ holds.  
\end{proof}

Let us also mention the following  construction from \cite{Pl02}.

\begin{theorem}
Under \CH\ there is   a nonmetrizable {\em convex}
compact subset $K$ of $\Sigma([0,1]^{\omega_1})$ carrying a strictly positive measure.
\end{theorem}

Such a space $K$ is defined as an affine  continuous image of some $P(L)$
(where $L$ is not Corson compact). Here one uses
the following observation: a compact space $L$ carries a strictly
positive measure if and only if  the space $P(L)$ carries such a measure, cf.\  
Marciszewski and Plebanek \cite{MP10}.

Coming back to first-countable spaces,
the following was proved in \cite{Pl00}.

\begin{theorem}\label{fc}
It is relatively consistent that  $P(K)$ is first-countable if (and only if) $K$ is 
first-countable for every compact space $K$.
\end{theorem} 

The assertion in \ref{fc} holds in the real random model and is a consequence  of
the fact that in this model $\omega_1$ is a caliber of Radon measures but, on the other hand, 
$2^\oo$ contains a subset of cardinality $\omega_1$ which has full outer measure.
However, the following problem posed by David H.\ Fremlin seems to be still open.

\begin{problem}
Does the assertion of \ref{fc} follows from $\MA(\omega_1)$?
\end{problem}

\newcommand{\azm}{\aleph_0{\rm -monolithic}}

A compact space $K$ is $\azm$ if every separable subspace of $K$ is metrizable.
Note that every Corson compact is $\azm$: if $K\sub \Sigma([0,1]^\Gamma)$ and
we take any sequence of $x_n\in K$ which is dense in $K_0\sub K$ then there is a countable set $\Gamma_0\sub\Gamma$
such that $x_n(\gamma)=0$ for every $n$ and every $\gamma\in\Gamma\sm\Gamma_0$.
Consequently, $K_0$ embeds into $[0,1]^{\Gamma_0}$.

In \cite{Pl21} we investigated compact spaces $K$ for which $P(K)$ is $\azm$. 

\begin{lemma}\label{mon:}
For any $K$ in \MS, the space $P(K)$ is $\azm$ if and only if
the support of every $\mu\in P(K)$ is metrizable.
\end{lemma}

\begin{proof}
Note first that the condition is sufficient: for any sequence of measures $\mu_n\in P(K)$ we let
$S$ be the support of $\mu=\sum_n 2^{-n}\mu_n$. Then $\mu_n\in P(S)\sub P(K)$ and $P(S)$
is metrizable since $S$ is such a space.

Consider now some $\mu\in P(K)$ and its support $S$. Since  $\mu$ has countable type, 
$P(S)$ is separable  by Corollary \ref{corta}. 
Hence $P(S)$ is a separable subspace of a monolithic space so $P(S)$ is metrizable
and thus $S$ is metrizable.
\end{proof}

\begin{corollary}\label{mon:0}
Assuming $\MA(\oo)$,   the space $P(K)$ is $\azm$ if and only if
the support of every $\mu\in P(K)$ is metrizable.
\end{corollary}

\begin{proof}
In view of Lemma \ref{mon:} it is enough to observe that if $K$ is not in \MS\ then
$P(K)$ cannot be $\azm$. This follows from Theorem \ref{fremlin} given in the next section.
\end{proof}

\begin{theorem}\label{mon:1} 
Under \MA$(\omega_1)$,  the following are equivalent for a compact space $K$

\begin{enumerate}[(i)]
\item  $P(K)$ is $\aleph_0$-monolithic;
\item $K$ is $\aleph_0$-monolithic;
\end{enumerate}
\end{theorem}

Here $(i)\to (ii)$ is obvious since monolithicity is a hereditary property.
The reverse implication follows from a more general result due to
 Arkhangel'ski\u{\i} and \ Shapirovski\u{\i} \cite{AS87},  stating that any $\azm$ compact space satisfying
 the countable chain condition is, under $\MA(\oo)$, metrizable.

The implication $(ii)\to (i)$ of  Theorem \ref{mon:1} is not provable in ZFC and is dramatically
violated by   \CH, see Theorem \ref{underch}.
Whether  Corollary \ref{mon:0} can be proved in ZFC was a question asked 
 by Wies{\l}aw Kubi\'s in connection with \cite{KK12} and \cite{FKK13}.
More recently, the same problem  was posed  by Claudia Correa, 
motivated by  \cite{Co19}. The negative answer came in \cite{Pl21} and reads as follows.

\begin{theorem}
Under $\diamondsuit$,  there is a nonmetrizable  Corson compact space $K$ such that $P(K)$ is
$\aleph_0$-monolithic but $K$ supports a measure of type $\omega_1$.
\end{theorem}

\section{Mappings onto Tikhonov cubes}\label{mtc}

If $f:K\to [0,1]^\kappa$ is a continuous surjection then $K$ carries a measure of type $\ge\kappa$. Indeed, by Theorem \ref{bp:2}
there is $\mu\in P(K)$ such that $f[\mu]=\lambda_\kappa$; note that then $\mt(\mu)\ge \kappa$.
Richard Haydon was first to address the question whether this can be reversed, if the existence of $\mu\in P(K)$ of type $\kappa$
yields a continuous surjection onto $[0,1]^\kappa$.
He proved  in \cite{Ha77} that this is the case 
for every $\kappa$ with the property that
$\tau<\kappa$ implies $\tau^{\om}<\kappa$; 
 $\kappa=\con^+$ is the first uncountable cardinal with such a property.
Assuming the continuum hypothesis,  Haydon \cite{Ha78} proved (in particular)
that there is $K$ such that $|K|=|P(K)|=\con$ and $K$ carries a measure of uncountable type.
Of course, such a space $K$ cannot be mapped onto $[0,1]^\con$. 
The same features has Kunen's space from \cite{Ku81} which is perfectly normal; see also
D\v{z}amonja and Kunen \cite{DK93} for related constructions.
On the other hand, the following holds.

\begin{theorem}[Fremlin \cite{Fr97}] \label{fremlin}
Under $\MA(\kappa)$, if $K$ carries a measure $\mu$ of type $\ge\kappa$ then
$K$ can be continuously mapped into $[0,1]^\kappa$.
\end{theorem}

The author proved in  \cite{Pl97} Theorem \ref{plebanek} stated below
(assuming $\kappa\ge {\rm cf}(\kappa)\ge\om_{2}$).
This result was extended by Richard Haydon  to any
$\kappa\ge\om_{2}$. We discuss here only cardinal numbers
having  uncountable cofinality, cf.\ \cite{Pl02a} for the remaining case. 
In fact Haydon singled out an interesting lemma,
 described by Fremlin \cite{Fr00}, which is stated as Lemma \ref{hflemma} in a more general form,
enabling one to prove also a related result due to Talagrand (Theorem \ref{talagrand}).
Note that, unlike the previously mentioned results,  \ref{talagrand} requires no
additional set-theoretic assumptions.

A family of disjoint pairs
$((A_{\alpha}^{0}, A_{\alpha}^{1}))_{\alpha<\kappa}$ of sets  is said to be
independent if 
\[ A(\phi)=\bigcap_{\alpha\in I} A_{\alpha}^{\phi(\alpha)}\not=\emptyset\]
for every finite $I\sub\kappa$ and $\phi:I\to\{0,1\}$.
In the case when all the sets $A_{\alpha}^{i}$ are measurable with respect
to some measure $\mu$, we say that a family of pairs
$((A_{\alpha}^{0}, A_{\alpha}^{1}))_{\alpha<\kappa}$ is $\mu$--independent
if
we always have $\mu(A(\phi))>0 $
 ($\mu$--independence should not be confused with stochastic
 independence).

Independent families are the basic tool for defining continuous
surjections onto Tikhonov  cubes, as  is explained in the following lemma
(see \cite{Ha77}, Lemma 1.1).

\begin{lemma}\label{independent}
If $K$ is a compact space and $\kappa$ is a cardinal
number then the following are equivalent

\begin{enumerate}[(i)]
\item there is a continuous surjection from $K$ onto $\zok$;
\item  there is an independent sequence
$((F_{\alpha}^{0}, F_{\alpha}^{1}))_{\alpha<\kappa}$ of disjoint pairs, where
$F_{\alpha}^{0}$ and $F_{\alpha}^{1}$
are closed subsets of $K$ for every $\alpha<\kappa$.
\end{enumerate}
\end{lemma}

In the Cantor cube $\ttk$ we denote by $C_{\alpha}^{i}$, for $\alpha<\kappa$ and
$i=0,1$,  the one dimensional cylinders, i.e.
\[C_{\alpha}^{i}=\{x\in\ttk:\; x(\alpha)=i\}.\]

\begin{lemma}(Haydon-Fremlin)\label{hflemma}
Let $\kappa$ be any cardinal number
such that $\kappa\ge \om_{2}$. Suppose that
$((A_{\alpha}^{0}, A_{\alpha}^{1}))_{\alpha<\kappa}$ is a sequence
of pairs of measurable subsets of $\ttk$ with the following properties:

\begin{enumerate}[(i)]
\item $A_{\alpha}^{i}\sub C_{\alpha}^{i}$ for every
$\alpha<\kappa$ and $i=0,1$;

\item $\mk(A_{\alpha}^{0})+\mk(A_{\alpha}^{1})>\frac{1}{2}+r$
for every $\alpha<\kappa$;

\end{enumerate}

where $r$ is a constant such that $0\le r<\frac{1}{2}$. Then there are
a sequence of measurable sets $(Z_{\alpha})_{\alpha<\kappa}$ and
$X\in [\kappa]^{\kappa}$ such that

\begin{enumerate}[(a)]
\item $Z_{\alpha}\sub A_{\alpha}^{0}\cup A_{\alpha}^{1}$ and
$\mk (Z_{\alpha})>2r$ for every $\alpha<\kappa$

\item  for every $I\in [X]^{<\om}$,
if $\mk(\bigcap_{\alpha\in I}Z_{\alpha})>0$ then
$((A_{\alpha}^{0}, A_{\alpha}^{1}))_{\alpha\in I}$ is a $\mk$--independent
finite sequence.
\end{enumerate}
\end{lemma}

\begin{proof}
Recall that $(\ttk,\oplus)$ is a compact topological group
if we denote by $\oplus$ the coordinate-wise addition mod 2, and
$\mk$ is the Haar measure of that group.
Let $s_{\alpha}:\ttk\to\ttk$ be a mapping defined by
$s_{\alpha}(x)=x\oplus e_{\alpha}$, where $e_{\alpha}(\beta)=1$ iff
$\alpha=\beta$. Then $s_{\alpha}$ is a measure preserving homeomorphism of $\ttk$.

For every $\alpha<\kappa$ we can find zero  sets
$Z_{\alpha}^{i}\sub A_{\alpha}^{i}$, $i=0,1$, such that
$\mk(Z_{\alpha}^{0})+\mk(Z_{\alpha}^{1})>\frac{1}{2}+r$.
Put $H_{\alpha}=Z_{\alpha}^{0}\cup Z_{\alpha}^{1}$ and
$Z_{\alpha}=H_{\alpha}\cap s_{\alpha}[H_{\alpha}]$.
Then
\[\mk (Z_{\alpha})=\mk(H_{\alpha}\cap s_{\alpha}[H_{\alpha}])=%
2\mk(H_{\alpha})-\mk(H_{\alpha}\cup s_{\alpha}[H_{\alpha}])>2r.\]

Now the key point is that, by the definition of $s_{\alpha}$, the set $Z_{\alpha}$ 
is determined by coordinates in some countable set $J_{\alpha}\sub\kappa\sm\{\alpha\}$.
Since $\kappa\ge\omega_2$, by Hajnal's Free Set Theorem (see \cite[44.3]{Er84}), 
we can find a free set $X\sub \kappa $ of size $\kappa$, that is 
$\alpha\notin J_{\beta}$ whenever $\alpha, \beta\in X$.

For any finite $I\sub X$ and any function $\phi:I\to \{0,1\}$ we have
\[ \mk(\bigcap_{\alpha\in I}A_{\alpha}^{\phi(\alpha)})\ge\mk(\bigcap_{\alpha\in I}Z_{\alpha}\cap%
C_{\alpha}^{\phi(\alpha)})=\frac{1}{2^{|I|}}\mk (\bigcap_{\alpha\in I}%
Z_{\alpha}),\]
where the last equality comes from the fact that the sets $\bigcap_{\alpha\in I} Z_\alpha$
and $\bigcap_{\alpha\in I} C_\alpha^{\phi(\alpha}$ are determined by disjoint sets of coordinates
(and $\mk$ is a product measure).
\end{proof}

\begin{corollary}\label{cor_to_hf}
Let $\mu$ be a homogenous Radon measure of type $\kappa\ge\om_{2}$
on a compact space $K$.
Then for every $0< c< 1$ there are
a sequence $((F_{\alpha}^{0}, F_{\alpha}^{1}))_{\alpha<\kappa}$
of pairs of disjoint closed subsets of $K$ and a
sequence $(G_{\alpha})_{\alpha<\kappa}$ of Borel subsets of $K$
such that

\begin{enumerate} [(i)]
\item $G_{\alpha}\sub F_{\alpha}^{0}\cup F_{\alpha}^{1}$ and
$\mk (G_{\alpha})>c$ for every $\alpha<\kappa$;
\item for every $I\in [\kappa]^{<\om}$,
if $\mu(\bigcap_{\alpha\in I}G_{\alpha})>0$ then
$((F_{\alpha}^{0}, F_{\alpha}^{1}))_{\alpha\in I}$ is $\mu$--independent.
\end{enumerate}
\end{corollary}

\begin{proof} 
We only sketch the basic idea:
Since $\mu$ is a homogenous measure of type $\kappa$,
then there is an isomorphism $\theta: \fB_{\kappa}\to \fA$
between the measure algebra $\fA$ of $\mu$ and the measure algebra
$\fB_{\kappa}$ of $\mk$. 
For every $\alpha<\kappa$
we may find a Borel set $B_{\alpha}\sub K$ such that
$B_{\alpha}\kr =\theta (C_{\alpha}^{0}\kr)$. Next,
we find closed sets $F_{\alpha}^{0}$, $F_{\alpha}^{1}$ such that
$F_{\alpha}^{0}\sub B_{\alpha}$,
$F_{\alpha}^{1}\sub K\sm B_{\alpha}$, and
$\mu(F_{\alpha}^{0})+\mu(F_{\alpha}^{1})>\frac{1}{2}+c/2$,
for every $\alpha<\kappa$. Then we move to $\fB_\kappa$,  
apply Lemma \ref{hflemma} with $r=c/2$ and come back to $\fA$, choosing the required sets
in $K$ by regularity.
\end{proof}

We are now ready to show how this machinery  works.

\begin{theorem}  \label{plebanek}
If  $\kappa\ge\om_{2}$ is a caliber of Radon measures
then every compact space $K$ carrying a measure of type $\kappa$ can be continuously mapped onto $\zok$.
\end{theorem}

\begin{proof}
Let $\mu\in P(K)$ be a homogenous measure of type $\kappa\ge\om_{2}$.
In the notation of Corollary \ref{cor_to_hf}, where $c>0$,
we have $\mu(G_{\alpha})\ge c$ so
we can find $X\in [\kappa]^{\kappa}$, such that the family
$(G_{\alpha})_{\alpha\in X}$ is centered. 

Then
the family $((F_{\alpha}^{0}, F_{\alpha}^{1}))_{\alpha\in X}$ is
independent by Corollary 4.2 and we conclude  applying Lemma \ref{independent}.
\end{proof}

\begin{theorem} [Talagrand] \label{talagrand}
For every $\kappa\ge\om_{2}$, if a compact space $K$ carries a measure $\mu$ of type $\kappa$ then
$P(K)$ can be continuously mapped onto $\zok$.
\end{theorem}

\begin{proof} 
We fix $c$ with $1>c>1/2$ and,
keeping the notation form Corollary \ref{cor_to_hf},
for every $\alpha<\kappa$ and $i=0,1$ we put
\[ M_{\alpha}^{i}=\{\nu\in P(K):\; \nu(F_{\alpha}^{i})\ge c\}.\]
In this way for every $\alpha$ we have defined a disjoint pair
$(M_{\alpha}^{0}, M_{\alpha}^{1})$
of closed subsets of $P(K)$, so the proof will be complete
if we check that $((M_{\alpha}^{0}, M_{\alpha}^{1}))_{\alpha<\kappa}$
is an independent family.

Take any finite $I\sub \kappa$ and a function $\phi:I\to \{0,1\}$;
denote $H_{\alpha}=F_{\alpha}^{\phi(\alpha)}$ for simplicity.
We want to check that 
$\bigcap_{\alpha\in I}M_{\alpha}^{\phi(\alpha)}\not=\emptyset$;
in view of Theorem \ref{kelley} it will suffice  to check that the intersection number
of the family $\{H_{\alpha}:\; \alpha\in I\}$ is $\ge c$.

For a function $f=\sum_{\alpha\in I}n_{\alpha}\chi_{H_{\alpha}}$ and its essential supremum $\rm esup$
we have
\[ {\rm esup}(f)\ge\int_K f d\mu \ge c\sum_{\alpha\in I}n_{\alpha}.\]
Writing  $H=\{t\in K:\; f(t)={\rm esup}(f)\}$ we have $\mu(H)>0$ (since $f$ is integer-valued).
It follows that  for
$J=\{\alpha\in I:\; \mu(H\sm G_{\alpha})=0\}$ ,
we have $\mu(\bigcap_{\alpha\in J} G_{\alpha})>0$.

Now Corollary \ref{cor_to_hf} comes into play: it follows that
$\bigcap_{\alpha\in J} H_{\alpha}\not=\emptyset$ and,
taking any $t \in \bigcap_{\alpha\in J} H_{\alpha}$, we conclude that
\[   \Big\| \sum_{\alpha\in I}n_{\alpha}\chi_{H_{\alpha}}\Big\| \ge %
\sum_{\alpha\in I}n_{\alpha}\chi_{H_{\alpha}}(t )\ge %
\sum_{\alpha\in J}n_{\alpha}={\rm esup}(f)\ge c\sum_{\alpha\in I}n_{\alpha},\]
and  the proof is complete.
\end{proof}

\section{Tightness}\label{new}

The \textit{tightness} of a topological space $X$, denoted here by $\tau(X)$, is the least cardinal number
such that for every $A\sub X$ and $x\in\overline{A}$ there is a set $A_0\sub A$ with $|A_0|\le \tau(X)$
and such that $x\in\ol{A_0}$.  There is a convex analogue of tightness which can be discussed for instance
in dual Banach spaces. Given a Banach space $X$, consider the dual unit ball $B_{X^\ast}$ equipped with the $weak^\ast$ topology.
Then we write $\ct(B_{X^\ast})=\omega$ and say that  the ball  has convex countable tightness if for every $A\sub B_{X^\ast}$ and $x^\ast\in\ol{A}$ there is a countable
set $A_0\sub A$ such that $x^\ast$ is in the closure of the convex hull of $A_0$. 

Recall that a Banach space $X$ has property (C) of Corson if for every family $\mathcal{C}$ of convex closed subsets of $X$ we have $\bigcap\mathcal{C}\neq\emptyset$ provided that every countable subfamily of $\mathcal{C}$ has nonempty intersection. 
Pol \cite{Po80,Po82} proved that a Banach space $X$ has property (C) if and only if 
$B_{X^\ast}$ has convex countable tightness. He also raised a question if 
this is further equivalent to saying that $B_{X^\ast}$ has countable tightness.
Mart\'{\i}nez-Cervantes and  Poveda \cite{MP23} have recently proved that

\begin{theorem}\label{gonzalo}
Assuming Proper Forcing Axiom, for every Banach space $X$ the space
$(B_{X^\ast}, weak^\ast)$ has countable tightness 
if and only if it has convex countable tightness.
\end{theorem}

Coming back to $P(K)$ spaces, we shall discuss the present status of 
 the following two questions.

\begin{problem}\label{mainproblem}

\begin{enumerate}[(A)]
\item Suppose that $P(K)$ has   countable tightness. Does this imply that $K\in$\MS?
\item Does convex countable tightness of $P(K)$ always  implies its countable tightness?
\end{enumerate}
\end{problem}

There are several reasons why such questions are  interesting and delicate. 
First note that by Theorem \ref{talagrand},  if $K$ carries a measure of type $\ge\omega_2$ then
$P(K)$ maps onto $[0,1]^{\omega_2}$ so 
\[ \tau(P(K)) \ge \tau ([0,1]^{\omega_2})=\omega_2,\]
 because
tightness is not increased by continuous surjection of compact spaces.
In other  words, higher analogues of \ref{mainproblem}(A) have a positive answer.

Secondly, under $\MA(\omega_1)$,  by Theorem \ref{fremlin},
 if a compact space $K$ admits a measure of uncountable type then $K$ can be continuously mapped onto
$[0,1]^{\omega_1}$, so in particular $\tau(P(K))\ge \tau(K)\ge \omega_1$. 
 It follows that Problem \ref{mainproblem}(A)
has a positive solution under  $\MA(\omega_1)$. The same effect is on Problem \ref{mainproblem}(B):
if $K\in$\MS\  then the results holds by Theorem \ref{crn} below; if $K\notin$ \MS\ then
the existence of a continuous surjection $g: K\to [0,1]^{\omega_1}$ gives
an isometric embedding $f\mapsto f\circ g$ of $C([0,1]^{\omega_1})$ into $C(K)$. It follows that
$C(K)$ does not have property (C) so $\ct(P(K))>\omega$.

There is a close connection between Problem \ref{mainproblem} (A) and (B), given by  the 
following result from \cite[Theorem 3.2]{FPCRN}.

\begin{theorem}\label{crn}
For every $K$ in \MS\ the following are equivalent
\begin{enumerate}[(i)]
\item $P(K)$ has convex countable tightness;
\item $P(K)$ has countable tightness.
\end{enumerate}	
\end{theorem}

Consequently, if a stronger version of Problem \ref{mainproblem}(A),
{\em does $\ct(P(K))=\omega$ imply $K\in$\MS},
has a positive solution then
countable tightness of $P(K)$ is indeed equivalent to $C(K)$ having property (C).

The following was proved by  Plebanek and Sobota \cite[Theorem 5.6]{PS15}.

\begin{theorem}\label{sobota}
If $P(K\times K)$ has convex countable tightness then $K\in$\MS.
\end{theorem}

Combining Theorem \ref{sobota}  with Theorem \ref{crn} we conclude that 

\begin{corollary}
For every $K$, $\tau(P(K\times K))=\omega$ if and only if $\ct(P(K\times K)) =\omega$.
\end{corollary}	

The next conclusion follows from Theorem \ref{sobota} and Theorem \ref{krupski}.
	
\begin{corollary}
	For every compact space $K$, either $P(K\times K)$ contains a ${G}_\delta$ point or $P(K\times K)$ has uncountable tightness.
\end{corollary}

Let us remark that  $K\in$\MS\ is equivalent to $K\times K\in$\MS.
Theorem \ref{sobota} solved partially Problem \ref{mainproblem}(B)
but the main argument relied heavily on the structure of product spaces.
  Avil\'es,   Martínez-Cervantes,   Rodr\'{\i}guez and   Rueda Zoca \cite[Corollary 7.8]{AMRR} proved
that Theorem \ref{sobota} is in fact a  particular case of their result on injective
tensor product of Banach spaces having property (C).

Note  that $\tau(P(K))\le \tau (P(K\times K))$ so Theorem \ref{sobota}  raises the question if the sharp inequality here is possible.
This was settled by a delicate construction due to Koszmider and Silber \cite{KS24} who showed that
it is relatively consistent that Problem \ref{mainproblem}(A) has a negative solution.

\begin{theorem}\label{kosz_sil}
Under $\diamondsuit$, there is a compact space $K$ carrying a measure of uncountable type and such that
$\tau(P(K))=\omega$.		
\end{theorem}

Here are several other questions regarding  sequential properties that, to our knowledge, are open:

\begin{problem}\label{kosz}
 Is it relatively consistent that	
\begin{enumerate}	
\item there is $K\notin$\MS\ such that $P(K)$ is hereditary separable?
\item
 there is $K\notin$\MS\ such that $P(K)$ is sequential?
\item  $\tau(P(K))=\omega$ whenever $\tau(K)=\omega$?
\item $P(K)$ is sequentially compact for every sequentially compact $K$?
\end{enumerate}
\end{problem}

The first two question above are repeated from
  \cite{KS24} while the last one arose in a corresponcence with
Niels Laustsen and Antonio  Acuaviva Huertos.  
In connection with \ref{kosz}(1), note that  \stevo\ \cite{To17} proved, under
$\diamondsuit$, that  there are nonmetrizable compacta $K$ such that
$P(K)$ is hereditary separable in all finite powers.

A compact space $K$ is said to be \textit{Rosenthal compact} if $K$ embeds into $B_1(X)$, the space of Baire-one functions on a Polish space $X$ equipped with the topology of pointwise convergence. The class of Rosenthal compacta is stable under taking countable product and, by
a result of Godefroy \cite{Go80},   if $K$ is Rosenthal compact, then so is $P(K)$. Moreover, 
Rosenthal compacta are Fr\'echet-Urysohn spaces (see Bourgain, Fremlin and Talagrand \cite{BFT}), hence they have countable tightness. This, together with Theorem \ref{sobota}, implies the following

\begin{corollary}
	If $K$ is Rosenthal compact, then every $\mu\in P(K)$ has countable type.
\end{corollary}

This fact,  announced by J.\  Bourgain in his PhD thesis, was proved by Todor\v{c}evi\'c \cite{To99} basing on
properties of Rosenthal compacta in forcing extensions; see also Marciszewski and Plebanek \cite{MP12}.  
However, the following seems to be open; see 
 \cite{MP12} for  a partial positive solution.

\begin{problem}
Let $K$ be Rosenthal compact; is every $\mu\in P(K)$ countably determined?	
	\end{problem}

\section{Converging sequences of measures}\label{sp}

Given $K$ and $\mu\in P(K)$, a sequence $(x_n)_n$ in $K$ is $\mu$-uniformly distributed if
\[ \frac{1}{n}\sum_{k=1}^ n \delta_{x_k}\to \mu \mbox{ that is } \lim_n \frac{1}{n}\sum_{k=1}^n f(x_k)=\int_K f\; {\rm d}\mu\mbox{ for every } f\in C(K).\] 

By Lemma \ref{bp:4}, the set ${\rm conv}\Delta_K$ of probability measures supported by finite sets is dense in $P(K)$.
Niederreiter \cite{Ni72} proved that for a given  $\mu\in P(K)$, the existence of $\mu$-uniformly distributed sequence is
equivalent to the fact that $\mu$ is a limit of a converging sequence in ${\rm conv}\Delta_K$ (see also \cite[491D]{Fr4}).
It follows that if, for instance, $\mu\in P(K)$ has a countable local base, then $\mu$ admits   $\mu$-uniformly distributed sequence.
In fact,  Mercourakis \cite{Me96} proved that such a property has every measure which is countably determined.
Uniformly distributed sequences may be found in various settings, see
for instance Fremlin \cite[491F]{Fr4} for a general result on product measures and \cite[491H]{Fr4}
stating that  the Haar measure of every  compact separable group has a uniformly distributed sequence.
In particular,  the measure $\lambda_\con$ on $2^\con$
admits such a sequence (see \cite[491G]{Fr4}). The following theorem due to
Fremlin \cite[491Q]{Fr4} is a
far reaching generalization of that classical result.

\begin{theorem} \label{fr_ud}
Every $\mu\in P(2^\con)$ has a uniformly distributed sequence.
\end{theorem}

The first step towards \ref{fr_ud} was done by Losert \cite{Lo78} who proved the results with $\omega_1$ replacing $\con$.
Then \ref{fr_ud} was proved in \cite{FP96} but assuming $\MA$.
Finally, Fremlin found a completely different argument in ZFC and we outline it here.

Recall that the asymptotic density $d(A)$ of a set $A\sub\omega$ is defined as
\[d(A)=\lim_{n\ge 1} |A\cap n|/n\] 
whenever the limit exists (here $n=\{0,1,\ldots, n-1\}$). The upper asymptotic density
$\overline{d}$ is defined accordingly, for every $A\sub\omega$,  by replacing $\lim$ with $\limsup$ in the formula above.
 
Write $\cD$ for the family of those $A$ for which $d(A)$ is defined and
$\cZ=\{A\sub\omega: d(A)=0\}$ for the density zero ideal. We consider the quotient structure
$\cD/\cZ$ denoting by $A\kr$ the equivalence class of the set $A$.
Note that we may as well treat  $d$ as a function on $\cD/\cZ$.

Although $\cD/\cZ$ is not a Boolean algebra,  Fremlin  \cite[491P]{Fr4} proved an interesting
theorem on embedding of measure algebras into it (the formulation
given here is slightly more general).

\begin{theorem}\label{david}
Let $\fA$ be any Boolean algebra of cardinality $\le\con$ and let $\mu\in P(\fA)$.
Then there is a Boolean homomorphism $\vf:\fA\to \cD/\cZ$ such that
$\mu(a)=d(\vf(a))$ for every $a\in\fA$.
\end{theorem} 

\begin{proof}
Observe first that we only have to prove the assertion for
$\lambda_\con \in P(\fB_\con)$. Indeed, for any $\mu\in P(\fA)$ we can suitably embed $\fA$ into the measure algebra
of the corresponding  Radon measure defined on $Bor(\ult(\fA))$; in turn that Radon measure is of type $\le\con$
so, by the Maharam theorem,  its measure algebra embeds into $\fB_\con$ by the measure preserving homomorphism.

We have already mentioned that the product measure $\lambda=\lambda_\con$ admits a uniformly distributed sequence $(x_n)_n$.
This immediately gives a Boolean homomorphism $\psi$ from $\clop(2^\con)$ into the power set of $\omega$ 
defined by
\[ \psi(V)=\{n\in\omega: x_n\in V\} \mbox{ for } V\in\clop(2^\con).\]
Note that, by uniform distribution, we actually have $\psi(V)\in\cD$ and
$d(\psi(V))=\lambda(V)$ for every clopen set $V$.

Writing $V\kr$ for the element of $\fB_\con$ corresponding to $V\in\clop(2^\con)$ 
and  $\fC=\{V\kr: V\in\clop(2^\con)\}$ we
have  a Boolean homomorphism 
\[ \vf_0: \fC\to \cD/\cZ \mbox{  defined by } \vf_0(V\kr)=(\psi(V))\kr.\]

We now  equip $\fB_\con$ with the metric $\rho(\cdot,\cdot)=\lambda(\cdot\bigtriangleup \cdot)$ 
 and $\cD/\cZ$ with the analogous metric $\sigma(\cdot ,\cdot)=\overline{d}(\cdot \bigtriangleup \cdot)$.  Then we can treat
 $\vf_0$ as an isometric embedding of the dense subspace $\fC$ of $\fB_\con$
 into $\cD/\cZ$. The main technical point is that  $\cD/\cZ$ is a complete metric space
 so $\vf_0$ extends uniquely to an isometry $\vf:\fB_\con\to\cD/\cZ$; see
 \cite[491I(b)]{Fr4} for details. Then it is not difficult to check 
 such an isometry is the desired  homomorphism --- note the Boolean operations are continuous.
 \end{proof}
 
It is worth recalling here that the problem if the measure algebra
of the Lebesgue measure can be embedded  into the power set of $\omega$ mod the ideal
of finite sets  is undecidable within the usual set theory, see Dow and Hart \cite{DH00}.

\begin{proof}(of Theorem \ref{fr_ud})
Take $\mu\in P(2^\con)$ and a homomorphism $\vf$ from its measure algebra $\fA$
 into $\cD/\cZ$ as in Theorem \ref{david}. For $\alpha<\con$ and $i\in\{0,1\}$ write
 $C_\alpha=\{x\in 2^\con: x(\alpha)=1\}$ and choose
 $N_{\alpha}$ such that $\vf(C_\alpha\kr)=N_\alpha\kr$.
 
 Now we define $x_n\in 2^\con$ setting 
 \[x_n(\alpha)=1 \mbox{ if and only if } n\in N_{\alpha}.\]
 Then the points $x_n\in 2^\con$ form a $\mu$-uniformly distributed sequence. 
 To see this consider, for instance, the set $V=C_\alpha\sm C_\beta$; then
 \[ \frac{|\{n< k: x_n\in V\}|}{k}=\frac{|k\cap (N_{\alpha}\sm N_\beta)|}{k} \to  d(N_\alpha\kr \sm N_\beta\kr)=\]
 \[= d(\vf(C_\alpha\kr)\sm \vf(C_\beta\kr))=d(\vf(C_\alpha\kr\sm C_\beta\kr))=d(\vf(V\kr))=\mu(V).\]
By an analogous argument  we can check that  the same holds for any basic clopen set $V$ and so for any $V\in\clop(2^\con)$.
 \end{proof}

Recall that a space $K$ is dyadic if it is a continuous image of $2^\kappa$ for some $\kappa$. 

\begin{corollary}
For every separable dyadic space $K$,  every $\mu\in P(K)$	has a uniformly distributed sequence.
\end{corollary}	

\begin{proof}
It is known that if a dyadic space $K$ is separable then there is a continuous surjection $f:2^\con \to K$.
Given $\mu\in P(K)$, there is $\nu\in P(2^\con)$ such that $f[\nu]=\mu$ (see Theorem \ref{bp:2}). Now
it is routine to check that if $(x_n)_n$ is a sequence in $2^\con$ which is $\nu$-uniformly distributed then
the sequence of $f(x_n)\in K$ is $\mu$-uniformly distributed.	
\end{proof}	

\newcommand{\seq}{\protect{\rm Seq}}

Writing $\seq^0(K)={\rm conv}\Delta_K$ and $\seq^{\alpha}$ for all the measures in $P(K)$
that are limits of converging sequences from $\bigcup_{\beta<\alpha}  \seq^\beta(K)$, the space 
\[ \seq(K)=\bigcup_{\alpha<\omega_1}\seq^\alpha(K)\]
is the sequential closure of all finitely supported probabilities.
We have briefly discussed spaces for which $P(K)=\seq^1(K)$; however
the structure of such a sequential closure may bo more complicated.
The following result comes from \cite{APR13}.

\begin{theorem}\label{apr|}
Under \CH\  there is a compact space $K$ such that $\seq^1(K)\neq \seq(K)=P(K)$.
\end{theorem}

Borodulin-Nadzieja and Selim \cite{BNS} obtained the following partial generalization of that result; however
in their construction $\seq(K)\neq P(K)$.

\begin{theorem} \label{bns}
For every $\alpha<\omega_1$ there is a compact space $K$ such that
\[ \seq^{\alpha+1}(K)=\seq^\alpha(K)\neq\bigcup_{\beta<\alpha} \seq^\beta(K).\]
\end{theorem}

It seems that no more is known on the hierarchy of $\seq(K)$; in particular, the following
problems are open.

\begin{problem}
\begin{enumerate}[(A)]
\item Is it relatively consistent that $\seq(K)=\seq^1(K)$ whenever $K$ is such a compact space that
$\seq(K)=\seq^n(K)$ for some natural number $n$?
\item Is there a space $K$ satisfying the assertion of \ref{bns} for which $\seq(K)=P(K)$? 
\end{enumerate}
\end{problem}

\section{Grothendieck spaces and the Efimov problem}\label{groef}

A Banach space $X$ is a Grothendieck space if every $weak^\ast$ converging sequence 
in $X^\ast$ converges weakly. Gonzalez and Kania \cite{GK21} offer an extensive
survey on various aspects of the class of Grothendieck Banach spaces. 

Recall that  $weak$ convergence of sequence of signed measures $\mu_n\in C(K)^\ast$ 
is characterized by the condition that the limit $\lim_n \mu_n(B)$ exists for
every $B\in Bor(K)$ (see e.g.\ Diestel \cite[Theorem 11]{Di84}). Therefore,  a Banach space of the form $C(K)$ is
Grothendieck if the convergence of a sequence of signed measures on continuous functions
implies its convergence on all Borel sets.
This implies that $C(K)$ is not Grothendieck whenever $K$ contains a sequence of distinct $x_n\in K$ converging 
to some $x\in K$:
in such a case $\delta_{x_n}\to\delta_x$ but $\delta_{x_n}(\{x\})=0$ infinitely often.

The Banach space $C(\beta\omega)\equiv \ell_\infty$ is a classical example of  a Grothendieck space. This
fact was generalized in various direction, in particular for the Stone spaces of Boolean algebras
with some weak versions of sequential completeness, see Haydon \cite{Ha81},  Koszmider and Shelah \cite{KS13} 
for concrete results and further references. Moreover,  the authors of \cite{KS13}
named the positive version of the Grothendieck property which amounts to saying that every convergent sequence
of measures in $P(K)$ is convergent on every Borel set in $K$.

Cembranos \cite{Ce84} and Freniche \cite{Fre84} proved   (abstract) results  implying  that
no space of the form $C(K\times L)$ is Grothendieck whenever $K$ and $L$ are infinite compact spaces.
 This result has been recently reproved in \cite{KSZ} by constructing in a given product space $K\times L$
 a nontrivial sequence of purely atomic signed measures converging to 0 in the $weak^\ast$ topology but not converging
 weakly. 
 We present  below a variation on the subject stating that $P(K\times L)$
 is never positively Grothendieck for such product spaces.

\begin{theorem}\label{add:3}
If $K$ and $L$ are infinite compact spaces then there is a sequence of $\mu_n\in P(K\times L)$
converging to some $\mu\in P(K\times L)$ having the support $S$  such that $\mu_n(S)=0$
for every $n$.
\end{theorem}

\begin{proof}
It is well-known that if $K$ is scattered and infinite then $K$ contains a nontrivial converging sequence;
this quickly yields the assertion. Therefore we can assume that both $K$ and $L$ are not scattered
and hence we can fix  continuous surjections $f:K\to [0,1]$ and $f':L\to [0,1]$.

We first consider a certain sequence of measures on $[0,1]^2$. Write $\lambda$ for the Lebesgue measure on
$[0,1]$ and $\lambda_2$ for the product measure. For any $n\ge 0$ and $1\le k\le 2^n$ denote
$I(n,k)=[(k-1)/2^n, k/2^n]$ and
\[ D_n=\bigcup_{k\le 2^n} I(n,k)\times I (n,k).\]
Then $\lambda_2(D_n)=1/2^n$; we let $\nu_n$ be the normalized restriction of $\lambda_2$ to the set $D_n$.
Note that every $\nu_n$ has $\lambda$ as the marginal distributions, that is
\[ \nu_n(B\times [0,1])=\nu_n([0,1]\times B)=\lambda(B)\]
 for every $B\in Bor[0,1]$.

Observe that the measures $\nu_n$ can be uniformly approximated on rectangles:
If $B,B'\sub [0,1]$ are  Borel sets  and for $\eps>0$ we  find  $\eps$-approximations 
$F$ and $F'$, that is $\lambda(B\btu F),\lambda(B'\btu F')<\eps$
 then 
\[\nu_n(B\times B')\le \nu_n(F\times F') +\nu_n((B\btu F)\times [0,1])+\nu_n([0,1]\times (B'\btu F'))\le \] 
\[ \le \nu_n(F\times F')+\lambda(B\btu F)+\lambda(B'\btu F')\le \nu_n(F\times F')+2\eps.\]

It is clear that every  $\nu_n$ vanishes on the diagonal and it is easy to see that $\nu_n$ converge 
to the Lebesgue measure put on the diagonal. Moreover, by the uniform approximation mentioned above we have
\medskip

\noindent{\sc Claim 1.}
 $\lim_n \nu_n (B\times B')=\lambda(B\cap B')$ for every $B,B'\in  Bor[0,1]$.
\medskip

Now we consider the surjection 
\[ g: K\times L\to [0,1]^2,  \quad g(x,y)=(f(x), f'(y)) \mbox{ for } (x,y)\in K\times L,\] and
take  measures $\mu_K\in P(K)$, $\mu_L\in P(L)$ such that $f[\mu_K]=\lambda=f'[\mu_L]$. By Remark \ref{add:2}
we can assume that the $\sigma$-algebras 
\[ \Sigma=\{f^{-1}[B]: B\in Bor[0,1]\} \mbox{ and }\Sigma'=\{f'^{-1}[B]: B\in Bor[0,1]\}\]
 are $\bigtriangleup$-dense in $Bor(K)$ and $Bor(L)$, respectively.
 Then we have for every $n$ the measure $\mu_n\in P(K\times L)$ defined analogously to the shape of $\nu_n$ i.e.\
 $\mu_n$ is the normalized restriction of $\mu_K\otimes\mu_L$ to the set $g^{-1}[D_n]$.
 We conclude directly from Claim 1 that
\medskip

\noindent {\sc Claim 2.} The limit $\lim_n\mu_n(A\times A')$ exists for every $A\in\Sigma$ and $A'\in\Sigma'$.
\medskip

But then, using the density of $\Sigma$ and $\Sigma$ and  the  argument on uniform approximation,
we extend Claim 2 to saying that $\lim_n\mu_n(A\times A')$ exists for any Borel rectangle $A\times A'$.
Note at this point that every continuous function $h$ on $K\times L$ can be uniformly approximated
by a linear combination of characteristic functions of Borel rectangles; hence
the integrals $\int_{K\times L} h\; {\rm d}\mu_n $ converge for any $h\in C(K\times L)$.
By compactness $\mu_n\to \mu\in P(K\times K)$; then $\mu$ is concentrated on
the closed set $G=\bigcap_n g^{-1}[D_n]$ while $\mu_n(G)=0$ because
every $\nu_n$ vanishes on the diagonal in $[0,1]^2$.
\end{proof}

The Efimov problem is a question whether every infinite compact space either contains a nontrivial converging sequence
or a copy of $\beta\omega$, see Hart \cite{Ha07} for a survey on the subject.
Recall that a number of counterexamples have been found in various models of set theory --- those are usually called 
{\em Efimov spaces} ---
but it is a big open problem if there are such examples in ZFC. Recall also
that a compact space contains $\beta\omega$ if and only if it can be continuously mapped onto $[0,1]^\con$.

\begin{theorem} [Talagrand \cite{Ta80}] \label{ta_efimov}
Under \CH\ there is a Efimov space $K$ such that $C(K)$ is Grothendieck.
\end{theorem}

Talagrand's space $K$ is strongly Efimov: not only $K$ contains no converging sequences but
$P(K)$ contains no `nontrivial' converging sequence. 
Note that one can weaken the original Efimov problem and ask if every infinite compact space $K$ either
contains a nontrivial converging sequence
or it can be continuously mapped onto $[0,1]^{\omega_1}$.
In view of Theorem \ref{fremlin},  this could be further relaxed to

\begin{problem}\label{ef}
	Does every infinite compact space $K$ either contains a nontrivial
	converging sequence or carries a measure of uncountable type (i.e.\ $K\notin$\MS)?
	\end{problem}
	
This, however, was again  refuted under \CH, see D\v{z}amonja and Plebanek \cite{DP07}.
Later Dow and Pichardo-Mendoza \cite{DM09} proved the following finer result.

\begin{theorem}\label{dm}
Under \CH\ there is a zero-dimensional  Efimov space $K$ without isolated points which is a limit of inverse system of length $\omega_1$ consisting 
of simple extensions of metric compacta.
\end{theorem}

The inverse system of spaces $K_\alpha$ is built on simple extensions if
for every $\alpha$ the bonding map $K_{\alpha+1}\to K_\alpha$ splits one point of $K_\alpha$ into two
and is injective otherwise. If the compacta in question are zero-dimensional then, in the Boolean language,
the resulting space $K$ is the Stone space of a minimally generated Boolean algebra, 
see Borodulin-Nadzieja \cite{PBN07} for details, where the following was proved.

\begin{theorem}\label{pbn}
If $K$ is a limit of inverse system of length $\omega_1$ consisting 
of simple extensions of metric compacta then every nonatomic measure $\mu\in P(K)$
is strongly countably determined.
\end{theorem}

We can now contemplate the unusual properties of $K$ from Theorem \ref{dm}:

\begin{enumerate}
\item No $x\in K$ is $G_\delta$ and therefore  no $\delta_x\in P(K)$ is $G_\delta$; however, every
nonatomic $\mu\in P(K)$ is $G_\delta$ (by Proposition \ref{fact:Pol}).	
\item There are no nontrivial converging sequences in $K$ but there are plenty `nontrivial' converging sequences in $P(K)$ since
 every nonatomic $\mu\in P(K)$ admits a uniformly distributed sequence (see section \ref{sp}).	
\item Every $\mu\in P(K)$ is countably determined but $K$ is not sequentially compact.
\end{enumerate}	

Those remarks were already noted in \cite{DP07} but in connection with
the usual Fedorchuk-like construction which requires $\diamondsuit$.
Banakh and Gabriyelyan \cite{BG24} considerably refined 
that  $\diamondsuit$-construction to obtain several delicate properties of the resulting space $P(K)$;
in particular, for their space $K$, the subspace of $P(K)$ of all nonatomic measures is sequentially compact.
Other refinements can be found in Sobota and Zdomskyy \cite{SZ19} and \cite{SZ23}. 
Further discussion on connection of the Efimov problem and properties
of Grothendieck spaces can be found in \cite{KS13}.

\section{Some remarks}

One can wonder if, for a given space $K$, the space $P(K)$ is homeomorphic to some familiar
object. Recall that if $K$ is an infinite metrizable compactum then $P(K)$ is
homeomorphic to the Hilbert cube $[0,1]^\omega$; see
Fedorchuk \cite[section 3]{Fe91} and references therein. Ditor and Haydon \cite{DH76} proved
that $P(2^\kappa)$ is homeomorphic to $[0,1]^\kappa$ if and only if $\kappa\le\oo$;  in fact
for $\kappa>\oo$ the space $P(2^\oo)$ is not an absolute retract.

 Not much is known about  pairs of compact spaces $K$ and $L$, for which the spaces
 $P(K)$ and $P(L)$ are homeomorphic. Note, however, that
Avil\'es and Kalenda \cite{AK09} developed an interesting technique enabling one to determine
possible homeomorphisms between $P(K)$ and $P(L)$ for $K$ and $L$ from some classes of scattered
compact spaces.

\end{document}